\documentclass[11pt]{article}

\usepackage[a4paper, total={17cm, 27cm}]{geometry}

\usepackage{newtxtext}
\usepackage[vvarbb]{newtxmath}

\usepackage[utf8]{inputenc}

\usepackage{bm}
\usepackage{todonotes}
\usepackage{subcaption}
\usepackage{tikz}
\usetikzlibrary{patterns}
\usepackage{pgfplots,pgfplotstable}
\pgfplotsset{compat=newest}
\definecolor{blau0} {RGB}{ 131 198 216} 
\definecolor{grau}  {RGB}{  0  84 159}
\definecolor{rot}   {RGB}{204   7  30}
\definecolor{blau2} {RGB}{  0  61 128}
\definecolor{grun2} {RGB}{  0 85   0}
\definecolor{rot2}  {RGB}{120   7  30}
\definecolor{gelb}  {RGB}{70 70 70}

\definecolor{blau}{RGB}{0 144 188}
\definecolor{newblue1}{RGB}{0 144 188}
\definecolor{newblue2}{RGB}{197 216 227}
\definecolor{newgreen1}{RGB}{0 144 118}
\definecolor{grun}{RGB}{255 137 0}
\definecolor{newgreen2}{RGB}{197 222 215}
\definecolor{neworange1}{RGB}{255 137 0}
\definecolor{neworange2}{RGB}{255 205 105}

\usetikzlibrary{arrows,positioning} 
\tikzset{
    >=stealth',
    punkt/.style={
           rectangle,
           rounded corners,
           draw=black, very thick,
           text width=6.5em,
           minimum height=4em,
           text centered},
    mellomboks/.style={
           rectangle,
           draw=black,
           text width=12em,
           minimum height=2em,
           text centered},
    mellomboks2/.style={
           rectangle,
           draw=black,
           text width=9em,
           minimum height=2em,
           text centered},
    startslutt/.style={
           circle,
           draw=black, thick,
           text width=3em,
           minimum height=2em,
           text centered},
    pil/.style={
           ->,
           thick,
           shorten <=2pt,
           shorten >=2pt,}
}

\usepackage[linesnumbered,ruled,vlined]{algorithm2e}
\newtheorem{remark}{Remark}
\DeclareMathOperator*{\argmin}{arg\,min}

\usepackage{authblk}

\title{An accelerated staggered scheme for phase-field modeling of brittle fracture}

\author[1]{Erlend Storvik\footnote{Corresponding author: erlend.storvik@uib.no}}
\author[1]{Jakub Wiktor Both}
\author[2]{Juan Michael Sargado}
\author[1]{Jan Martin Nordbotten}
\author[1]{Florin Adrian Radu}

\affil[1]{Department of Mathematics, University of Bergen, All\'egaten 44, 5007 Bergen, Norway}
\affil[2]{Danish Hydrocarbon Research and Technology Centre, Technological University of Denmark, Elektrovej Bygning 375, 2800 Kgs Lyngby, Denmark}

\date{}

\begin{document}

\maketitle

\begin{abstract}
There is currently an increasing interest in developing efficient solvers for phase-field modeling of brittle fracture. The
governing equations for this problem originate from a constrained minimization of a non-convex energy functional, and the most commonly used solver is a staggered solution scheme. This is known to be robust compared to the monolithic Newton method, however, the staggered scheme often requires many iterations to converge when cracks are evolving. The focus of our work is to accelerate the solver through a scheme that sequentially applies Anderson acceleration and over-relaxation, switching back and forth depending on the residual evolution, and thereby ensuring a decreasing tendency. The resulting scheme takes advantage of the complementary strengths of Anderson acceleration and over-relaxation to make a robust and accelerating method for this problem. The new method is applied as a post-processing technique to the increments of the solver, hence, the implementation can be done with minor modifications to already available software. Moreover, the cost of combining the two acceleration schemes is negligible. The robustness and efficiency of the method are demonstrated through numerical examples.\end{abstract}

\section{Introduction}
Mathematical modeling of brittle fracture propagation is an important and challenging topic in engineering sciences. The main difficulty arises in the transition between the distinct material properties in the fracture and the bulk domain. In this paper, we consider a variational phase-field model, as introduced by Bourdin, Francfort, and Marigo \cite{francfort, francfortnumerics}. A smooth indicator that marks the broken and unbroken parts of the material regularizes the sharp crack topology. This enables modeling of fractures without conforming meshes or path-tracking algorithms (as in XFEM \cite{moes}). However, fine meshes are needed to resolve the regularized region between the fracture and the bulk domain.


The system is modeled by minimizing its energy as a function of material displacement and the indicator function. This leads to a system of coupled, nonlinear equations which is challenging to solve. The most common technique, due to its robust nature, is the staggered scheme. This method decouples the system and sequentially updates the displacement and indicator variable by solving their respective subproblems. However, the convergence properties are at times very bad, and iterating to satisfactory precision can result in large numbers of iterations \cite{gerasimovline, Farrell}. The monolithic Newton method, on the other hand, does not show the same numerical robustness. Therefore, several attempts have been made to find a method that is both fast and robust. A monolithic, modified Newton method was proposed in \cite{wickmodifiednewton}, a monolithic quasi-Newton method of BFGS type was applied in \cite{kristensenquasinewton} and \cite{wu}, a monolithic line-search Newton method (dependent on the system energy) was applied in \cite{gerasimovline}, and the truncated nonsmooth Newton multigrid method was proposed in \cite{graser}. In \cite{brun}, the L-scheme \cite{pop,list} was applied in the context of an augmented Lagrangian solver, and a combination of an over-relaxed staggered scheme and the monolithic Newton method was applied in \cite{Farrell}. 

 In this paper, we propose a novel strategy to accelerate the classical staggered solution scheme solely utilizing two techniques for post-processing increments: Anderson acceleration and over-relaxation. In addition to accelerating the staggered scheme without sacrificing robustness, the new method allows the use of already available staggered scheme solvers with minor modifications to the implementation.

%

Anderson acceleration was first developed in \cite{anderson} for integral equations. Since then, it has seen many applications, including electronic structure computations \cite{fang} and flow in deformable porous media \cite{jakubanderson}. It is a multi-secant, quasi-Newton method that has been related to a preconditioned GMRES \cite{walker}. Moreover, the method post-processes the increments of the solver by approximating the inverse of the Jacobian of the system by reusing previous iterations. It can, therefore, easily be applied in combination with splitting techniques such as the staggered scheme while maintaining the decoupled nature of the scheme. 

In \cite{pollock}, the authors show theoretically that the Anderson acceleration improves the convergence rate of linearly convergent schemes, which is the case of the staggered scheme. However, as proved in \cite{kelley}, the convergence is only local. In the case of phase-field modeling of brittle fracture, this is a challenge; when fractures are initiating or propagating, the system state usually jumps drastically between consecutive loading steps. Therefore, a ``naive'' application of Anderson acceleration is not suitable for this application, as will be demonstrated in the numerical examples of this work. Recently, there has been an increasing interest in modified Anderson acceleration methods to overcome issues of local convergence. In \cite{AAsafetyguards} a safeguard, based on the residual norm of the problem, is applied to restart Anderson acceleration, and in \cite{restartAA} a periodically restarted Anderson acceleration is applied within a Richardson fixed-point iteration to accelerate the convergence of iterative solvers for large sparse linear systems. 

Relaxation was applied to the staggered scheme on a phase-field model of brittle fracture in \cite{Farrell}. It is a post-processing method that updates each iterate by relaxing (scaling) its increment. For the purpose of this work, over-relaxation (a scaling larger than one) is of particular interest. This is because the staggered scheme steadily moves towards the final configuration of each loading step, and over-relaxation might move the iterates further during each iteration, potentially accelerating the convergence. For the particular loading steps in which fractures are propagating, the gain can be quite substantial. There is, however, a drawback with over-relaxation: Near the solution of each loading step one might end up over- and undershooting the solution sequentially leading to poor performance.

The most important observation of this paper is the complementary strengths of these two acceleration techniques; Anderson acceleration accelerates close to the solution, while over-relaxation accelerates during loading steps with large jumps in the solution (e.g., during crack propagation). We propose an acceleration algorithm that switches between Anderson acceleration and over-relaxation during each loading-step ensuring convergence at an accelerated rate. This scheme is related to the one in \cite{Farrell} where the authors switch between over-relaxation and monolithic Newton. However, for the new acceleration scheme, proposed in this paper, both of the combined acceleration methods function as post-processes to the increments of the standard staggered scheme. In other words, the new acceleration method can be implemented with minor modifications to already available software. Moreover, switching between the two acceleration techniques does not change the sparsity of the underlying linear systems. The switch criterion is based on the history of the residual norms of the staggered solution steps.  

To summarize, the main contributions in this paper are:
\begin{itemize}
 \item Presentation of the difficulties encountered with the application of plain Anderson acceleration and over-relaxation applied to the staggered solution scheme for variational phase-field modeling of brittle fracture. 
 \item A new acceleration algorithm that exploits the complementary strengths of Anderson acceleration and over-relaxation, utilizing residual norm evolution as a rule for switching between the methods.
 \item The performance of the proposed acceleration scheme is demonstrated through thorough numerical examples including classical benchmark problems. 
\end{itemize}

The paper is structured as follows: The mathematical model and numerical discretization are presented in Section~\ref{sec:math}. Here, we introduce the energy functional which is subject to minimization together with the discretization. In Section~\ref{sec:acceleration}, the staggered scheme and the acceleration techniques are presented. Both Anderson acceleration and relaxation are described before the combined acceleration scheme is presented together with the inexact Newton modification. Section~\ref{sec:numerics} contains the numerical study of the accelerations applied to the staggered scheme. We test the staggered scheme both with and without the combinations of Anderson acceleration and relaxation. Moreover, the optimal depth of Anderson acceleration and the choice of relaxation parameter is discussed. Finally, some concluding remarks are made in Section~\ref{sec:conclusion}.

\section{Mathematical problem}\label{sec:math}
In this section, the mathematical problem that is considered throughout the paper is presented. An elastic medium, represented by the domain $\Omega\subset \mathbb{R}^d$ with $d= 2$ or $3$, is subject to loading through traction forces, $\bm t$, along $\Gamma_N$ and displacement, $\bm u_D$, along $\Gamma_D$ to the extent that it might break. Here, $\Gamma_N\cup\Gamma_D=\partial\Omega$ are subsets of the boundary of the domain, $\Omega$, and $\Gamma_D$ has nonzero measure. The state of the material is modeled by Griffith's criterion \cite{griffith}, with constant $G_c$, and a smooth indicator function (the phase-field variable) $\varphi:\Omega\rightarrow [0,1]$ describes the state of the damage to the material. The phase-field is defined to take the value $0$ whenever the material is unbroken, and $1$ when the material is broken, and a model parameter $\ell$ determines the width of the regularized zone where the phase-field transitions from $0$ to $1$.

\subsection{The energy of the system}
Following the work of \cite{francfort,bourdin}, we can express the total energy of the system as a sum of the medium's elastic energy, the surface energy dissipation associated with the broken parts of the material and external work related to traction. 
Now, let $\bm u$ denote the material displacement and define the total energy functional as 
\begin{equation}
 \label{eq:energy}
 \mathcal{E}(\bm u, \varphi):=\int_\Omega\mathcal{E}_c(\varphi)+ \mathcal{E}_m(\bm u,\varphi)  \; d{\bf x}-\int_{\Gamma_N} \bm t \cdot \bm u\; ds
\end{equation}
where 
\begin{equation}\label{eq:phasefield}
 \mathcal{E}_c(\varphi) := \frac{G_c}{2}\left(\frac{\varphi^2}{\ell}+\ell\nabla\varphi\cdot\nabla\varphi\right),
\end{equation}
and 
\begin{equation}
 \label{eq:mechanics}
 \mathcal{E}_m(\bm u, \varphi) := g(\varphi){\Psi}^+(\bm \varepsilon)+ {\Psi}^-(\bm\varepsilon) -\bm b\cdot \bm u.
\end{equation}
Here, we have applied the degradation function $$g(\varphi):=(1-\kappa)(1-\varphi)^2+\kappa,$$ where $\kappa$ is a ``small'' constant. Other choices have been proposed in \cite{sargado}. Moreover, the material is assumed to be homogeneous and isotropic, and the elastic strain energy functional \begin{equation}\label{eq:elasticenergy}
{\Psi}(\bm \varepsilon) :=\frac{1}{2}\bm\varepsilon \!:\!\mathbb{C}\!:\!\bm\varepsilon =\mu(\bm\varepsilon\!:\!\bm\varepsilon)+\frac{\lambda\mathrm{tr}\left(\bm\varepsilon\right)^2}{2},
\end{equation}
where $\bm \varepsilon = \frac{\nabla \bm u+\nabla \bm u^\top}{2}$ is the linearized elastic strain tensor and $\mu$ and $\lambda$ are the Lam\'e parameters, has been decomposed into ``tensile'', ${\Psi}^+$, and ``compressive'', ${\Psi}^-$, parts. The additive spectral decomposition
$$\Psi^\pm\left(\bm\varepsilon\right):=\mu(\bm\varepsilon_\pm\!:\!\bm\varepsilon_\pm)+\frac{\lambda\langle\mathrm{tr}\left(\bm\varepsilon\right)\rangle_\pm^2}{2},$$
proposed in \cite{miehesplit} has been employed.
Here, $\langle a \rangle_\pm :=\tfrac{1}{2}(a\pm|a|)$ and $\bm \varepsilon_\pm :=\sum_{i}\langle \varepsilon_i\rangle_\pm\bm n_i \otimes \bm n_i$ where $\{\varepsilon_i\}$ and $\{\bm n_i\}$ are the principal strain and principal strain directions, respectively. Additionally, the material is unable to heal, and the constraint $\partial_t \varphi \geq 0$ is applied accordingly.

\subsection{Time discretized, contiuous-in-space equations}
The loading procedure is discretized by the implicit Euler scheme, giving the non-healing constraint at loading step $n\geq1$: 
\begin{equation}\label{eq:nonhealing}
 \varphi^n(x)-\varphi^{n-1}(x) \geq 0 \quad \forall x\in\Omega .
\end{equation}
Now, we define the displacement solution space
$ \bm V^n = \left\{\bm v \in \left(H^1(\Omega)\right)^d\ \big\vert\ \bm  v|_{\Gamma_D}=\bm u_D^n\right\},$ the displacement test space $\bm V^0=\left\{\bm v \in \left(H^1(\Omega)\right)^d\ \big\vert\ \bm  v|_{\Gamma_D}=0\right\}$, and the phase-field solution and test space 
$Q = H^1(\Omega)$.
Then, the solution $(\bm u^n, \varphi^n)\in \bm V^n \times Q$ at loading step $n\geq1$ is given by 
\begin{equation}
 \label{eq:solution}
 (\bm u^n,\varphi^n):=\argmin_{\bm u, \varphi}\left\{ \mathcal{E}(\bm u, \varphi, \bm t^n)\; |\; \bm u \in\bm V^n,\varphi\in Q\right\}.
\end{equation}
Letting $\langle \cdot,\cdot\rangle_X$ denote the usual $L^2$ inner product over the domain $X$ and denoting $$\bm\sigma^\pm(\bm u ) := \dfrac{\partial{\Psi}^\pm\left(\bm\varepsilon(\bm u)\right)}{\partial\bm\varepsilon(\bm u)},$$ we find the variation of the energy \eqref{eq:energy} with respect to $\bm u$ and $\varphi$ respectively: 
\begin{eqnarray}\label{eq:variationalmechanics}
\mathcal{E}_{\delta \bm u}(\bm u,\varphi,\bm v) &=& \left\langle \left(g(\varphi)\bm\sigma^+(\bm u)+\bm\sigma^-(\bm u)\right), \bm \varepsilon(\bm v)\right\rangle_\Omega -\left\langle\bm b,\bm v\right\rangle_\Omega-\left\langle \bm t,\bm v\right\rangle_{\Gamma_N}  \\
\mathcal{E}_{\delta \varphi}(\bm u, \varphi, q) &=&  \left\langle g'(\varphi){\Psi}^+(\bm \varepsilon),q\right\rangle_\Omega + \frac{G_c}{\ell}\left(\left\langle \varphi,q\right\rangle_\Omega +\ell^2\left\langle\nabla \varphi,\nabla q\right\rangle_\Omega\right).
\end{eqnarray}
It is now easy to see that the solution to \eqref{eq:solution}, $\left(\bm u^n, \varphi^n\right)$, satisfies the system of equations
\begin{eqnarray}
 \mathcal{E}_{\delta \bm u}(\bm u^n,\varphi^n,\bm v) &=& 0\\
 \mathcal{E}_{\delta \varphi}(\bm u^n, \varphi^n, q) &=& 0
\end{eqnarray}
for all $\bm v \in \bm V^0$ and $q\in Q$.

The inequality \eqref{eq:nonhealing} still requires some special treatment, and in this paper, we follow the approach of \cite{miehesplit} and introduce a history variable; 
\begin{equation}\label{eq:history}
 \mathcal{H}^n := \max_{k\leq n}{\Psi}^+\left(\bm \varepsilon\left(\bm u^k\right) \right).
\end{equation}
A modified version of the variation with respect to $\varphi$ is defined by
\begin{equation}
 \widetilde{\mathcal{E}}_{\delta \varphi}(\bm u, \varphi, q) =  \left\langle g'(\varphi)\mathcal{H}^n,q\right\rangle_\Omega + \frac{G_c}{\ell}\left(\left\langle \varphi,q\right\rangle_\Omega +\ell^2\left\langle\nabla \varphi,\nabla q\right\rangle_\Omega\right),
\end{equation}
and the solution at loading step $n$ will be defined as the pair $\left(\bm u^n,\varphi^n\right)\in \bm V^n\times Q$ that satisfies 
\begin{eqnarray}\label{eq:solution1}
 \mathcal{E}_{\delta \bm u}(\bm u^n,\varphi^n,\bm v) &=& 0\\
 \widetilde{\mathcal{E}}_{\delta \varphi}(\bm u^n, \varphi^n, q) &=& 0\label{eq:solution2}
\end{eqnarray}
for all $\bm v\in \bm V^0$ and $q \in Q$.
\subsection{Spatial discretization}
To solve \eqref{eq:solution1}--\eqref{eq:solution2} we apply conforming linear finite elements \cite{francfortnumerics, gerasimovline, Farrell, sargado}, both for the phase-field variable and the displacement. Let $\mathcal{T}^h_\Omega=\left\{T_k\right\}_k$ be a decomposition of the domain $\Omega$ into simplices, $T_k$, and define, at loading step $n\geq1$, the spaces 
\begin{eqnarray}
 \bm V_h^n &=& \left\{\bm v_h \in (H^1(\Omega))^d \; \big\vert\; \bm v_h|_{T_k}\in (\mathcal{P}^1(T_k))^d\; \forall\; T_k \in \mathcal{T}_\Omega, \; \bm v_h|_{\Gamma_D}=\bm u_D^n \right\},\nonumber \\
 Q_h &=& \left\{q_h \in H^1(\Omega) \; \big\vert\;  q_h|_{T_k}\in \mathcal{P}^1(T_k)\; \forall\; T_k \in \mathcal{T}_\Omega \right\},\nonumber
\end{eqnarray}
and $\bm V_h^0$ accordingly with zero trace.
The system of equations to be solved is then: Find $(\bm u^n_h,\varphi^n_h)\in {\bm V^n_h \times Q_h}$ such that 
\begin{eqnarray}\label{eq:solutiondisc1}
 \mathcal{E}_{\delta \bm u}(\bm u^n_h,\varphi^n_h,\bm v_h) &=& 0\\
 \widetilde{\mathcal{E}}_{\delta \varphi}(\bm u^n_h, \varphi^n_h, q_h) &=& 0\label{eq:solutiondisc2}
\end{eqnarray}
for all $\bm v_h \in \bm V^0_h$ and  $q_h \in Q_h$. Following standard procedures, this naturally translates to the algebraic residual equations \begin{eqnarray}
                                                                                                                                                            \mathrm{Res}_{\bm u}\left(\bm u_h^n,\varphi_h^n\right)&=&0\\
                                                                                                                                                            \mathrm{Res}_{\varphi}\left(\bm u_h^n, \varphi_h^n\right)&=&0,
                                                                                                                                                               \end{eqnarray}
where $\mathrm{Res}_{\bm u}$ and $\mathrm{Res}_{\varphi}$ denote the algebraic residuals corresponding to \eqref{eq:solutiondisc1} and \eqref{eq:solutiondisc2}, respectively.

\section{Staggered scheme and acceleration}\label{sec:acceleration}
The discrete governing equations \eqref{eq:solutiondisc1}--\eqref{eq:solutiondisc2} are strongly nonlinear and coupled. In this paper, we apply the staggered scheme \cite{Farrell, gerasimov, sargado} to solve them, decoupling the equations. We let $i\geq 1$ be the iteration index and define the staggered scheme as: Given $\varphi^{n,i-1}_h\in Q_h$, find $(\bm u^{n,i}_h,\varphi^{n,i}_h)\in {\bm V^n_h \times Q_h}$ such that 
\begin{eqnarray}\label{eq:solutionalternate1}
 \mathcal{E}_{\delta \bm u}(\bm u^{n,i}_h,\varphi^{n,i-1}_h,\bm v_h) &=& 0\\
 \widetilde{\mathcal{E}}_{\delta \varphi}(\bm u^{n,i}_h, \varphi^{n,i}_h, q_h) &=& 0\label{eq:solutionalternate2}
\end{eqnarray}
for all $(\bm v_h,q_h)\in {\bm V^0_h \times Q_h}$ and $\varphi^{n,0}_h := \varphi^{n-1}_h$. The iterations are terminated when the following stopping criterions are reached:
\begin{eqnarray}
 \left\|\mathrm{Res}_{\bm u}\left(\bm u_h^{n,i},\varphi_h^{n,i}\right)\right\|_2&\leq& \mathrm{Tol}_\mathrm{Res,Abs},\label{eq:stopresabs} \\ \label{eq:relresidual}
  \frac{\left\|\mathrm{Res}_{\bm u}\left(\bm u_h^{n,i},\varphi_h^{n,i}\right)\right\|_2}{\left\|\mathrm{Res}_{\bm u}\left(\bm u_h^{n,1},\varphi_h^{n,0}\right)\right\|_2}&\leq& \mathrm{Tol}_\mathrm{Res,Rel}, \\
  \left\|\bm u_h^{n,i}-\bm u_h^{n,i-1}\right\|_{L^2(\Omega)}+\left\|\varphi_h^{n,i}-\varphi_h^{n,i-1}\right\|_{L^2(\Omega)} &\leq& \mathrm{Tol}_\mathrm{Inc,Abs},\\
    \frac{\left\|\bm u_h^{n,i}-\bm u_h^{n,i-1}\right\|_{L^2(\Omega)}}{\left\|\bm u_h^{n,1}\right\|_{L^2(\Omega)}}+\frac{\left\|\varphi_h^{n,i}-\varphi_h^{n,i-1}\right\|_{L^2(\Omega)}}{\left\|\varphi_h^{n,0}\right\|_{L^2(\Omega)}} &\leq& \mathrm{Tol}_\mathrm{Inc,Res},
\end{eqnarray}
for given tolerances $\mathrm{Tol}_\mathrm{Res,Abs}$, $\mathrm{Tol}_\mathrm{Res,Rel}$, $\mathrm{Tol}_\mathrm{Inc,Abs}$ and $\mathrm{Tol}_\mathrm{Inc,Rel}$.
Notice that controlling the residuals corresponding to the phase-field equation \eqref{eq:solutionalternate2} is redundant due to it being solved second in the staggered scheme by an exact linear solver. 

To solve the nonlinear equation \eqref{eq:solutionalternate1} we apply the Newton method with the relative stopping criterion 
\begin{equation}\label{eq:innernewton}\frac{\left\|\mathrm{Res}_{\bm u}\left(\bm u_h^{n,i,j},\varphi_h^{n,i-1}\right)\right\|_2}{\left\|\mathrm{Res}_{\bm u}\left(\bm u_h^{n,1},\varphi_h^{n,0}\right)\right\|_2}\leq\mathrm{Tol}_\mathrm{inner}.\end{equation}
Here, $j\geq 1$ is the iteration index for the Newton method and the initial guess is chosen as the previous staggered iteration $\bm u_h^{n,i,0}:=\bm u_h^{n,i-1}$.                                                                                                                                                                                                                                                                                                         


The staggered scheme \eqref{eq:solutionalternate1}--\eqref{eq:solutionalternate2} is closely related to the alternate minimization method (it differs in the application of the history variable \eqref{eq:history}) and is known to be a robust solution method \cite{bourdinAM}. However, it might require a large number of iterations to reach satisfactory tolerances \cite{gerasimovline, wu, Farrell}. We aim to accelerate this slow convergence and propose a combination of Anderson acceleration and over-relaxation. We note that the staggered solution scheme can be written as the fixed-point iteration
\begin{equation}\label{eq:fixpoint}
 {\bf x}^{n,i}_h := \mathcal{S}({\bf x}^{n,i-1}_h) = {\bf x}^{n,i-1}_h+\Delta\mathcal{S}({\bf x}_h^{n,i-1})
\end{equation}
where $\mathcal{S}$ is the staggered solution scheme operator, $\Delta\mathcal{S}$ is the increment of the staggered scheme and ${\bf x}_h^{n,i}$ is the vector $\begin{pmatrix}\bm u^{n,i}_h\\ \varphi^{n,i}_h\end{pmatrix}$.

Now, we present both the Anderson acceleration and the relaxed staggered scheme and describe their strengths and weaknesses. Then, taking advantage of the strengths of both schemes, a combined scheme is presented. 

\subsection{Anderson acceleration}\label{sec:aa}
Anderson acceleration is a multi-secant method that mimics the monolithic Newton method. The acceleration acts as a post-processing procedure that updates the current iterate by a linear combination of the $m$ previous iterates, according to their respective increments. The value of $m$ is free to be chosen and is known as the depth of the acceleration. Moreover, it can be chosen adaptively. At loading step $n$, the Anderson accelerated staggered scheme of depth $m$ reads:\\
\begin{algorithm}[H]
Given ${\bf x}^0$\;
\For{$i = 1,2,...$ until convergence}{
Set depth $m_i=\min\{m,i-1\}$\;
Define ${\bf F}^i:=\left[\Delta\mathcal{S}\left({\bf x}_h^{n,i-m_i-1}\right),...,\Delta\mathcal{S}\left({\bf x}_h^{n,i-1}\right)\right]$\;
Let ${\bm\alpha}^{i}=\left[\alpha^i_{0},...,\alpha^i_{m_i} \right]^\top\in \mathbb{R}^{m_i+1}$ be the minimizer of $\left\| {\bf F}^i\bm\alpha^{i}\right\|_2$ subject to $\sum_k \alpha_k^i =1$\;
Define the accelerated iterate ${\bf x}^i_h:=\sum_{k=0}^{m_i}\alpha_k\mathcal{S}\left({\bf x}_h^{n,k+i-m_i-1}\right)$
 }
 \caption{Anderson acceleration}
 \label{alg:AA}
\end{algorithm}

Algorithm~\ref{alg:AA} is independent of the underlying fixed-point iteration, but is presented for the application to the staggered scheme here. An important feature of Anderson acceleration is that it preserves the decoupled nature of the staggered scheme, hence, the subproblem solvers are unaffected by it.

It is demonstrated in the numerical section that Anderson acceleration improves the convergence when close to the solution. However, the acceleration might deteriorate otherwise. This is especially important to notice for brutal crack propagation, where it sometimes fails to converge at all.

\subsection{Over-relaxation}
Relaxation applied to each subproblem of the staggered solution scheme was described and applied in \cite{Farrell}. The method first calculates the increment $\Delta\bm u^{n,i-1}_h$ obtained by solving equation \eqref{eq:solutionalternate1},  before defining the updated iterate as $$\bm u_h^{n,i}:=\bm u^{n,i-1}_h +\omega \Delta \bm u^{n,i-1}_h,$$
where $\omega\in(0,2)$ is a parameter. This new iterate ${\bm u}^{n,i}_h$ is now passed on to equation~\eqref{eq:solutionalternate2} and the same procedure is executed for the phase-field resulting in the updated iterate $\varphi^{n,i}_h$. Following standard literature on iterative methods, we refer to the choice $\omega\in(1,2)$ as over-relaxation and $\omega\in(0,1)$ as under-relaxation. At the $n$-th loading step the relaxed staggered scheme reads:\\
\begin{algorithm}[H]
Given $\varphi^{n,0}_h$ and $\omega\in(0,2)$\;
\For{$i = 1,2,...$ until convergence}{
Find $ \hat{\bm u}^{n,i}_h\in \bm V_h^n$ satisfying $\mathcal{E}_{\delta \bm u}(\hat{\bm u}^{n,i}_h,\varphi^{n,i-1}_h,\bm v_h) = 0,\quad \forall \bm v_h\in V_h^0$\;
Define $\Delta \bm u^{n,i-1}_h:=\hat{\bm u}^{n,i}_h-\bm u^{n,i-1}_h$\;
Update the iterate $\bm u^{n,i}_h:=\bm u^{n,i-1}_h +\omega \Delta \bm u^{n,i-1}_h$\;
Find $\hat{\varphi}^{n,i}_h\in Q_h$ satisfying $\widetilde{\mathcal{E}}_{\delta \varphi}(\bm u^{n,i}_h,\hat{\varphi}^{n,i}_h,q_h) = 0,\quad \forall q_h\in Q_h$\;
Define $\Delta \varphi^{n,i-1}_h:=\hat{\varphi}^{n,i}_h-\varphi^{n,i-1}_h$\;
Update the iterate $\varphi^{n,i}_h:=\varphi^{n,i-1}_h +\omega \Delta \varphi^{n,i-1}_h$\;
}
\caption{Relaxed staggered scheme}
\label{alg:OR}
\end{algorithm}

Under-relaxation is robust when applied to the staggered scheme, however, it usually slows down the scheme. Over-relaxation, on the other hand, tends to accelerate the loading steps of the staggered solution scheme where cracks occur, while it might slow down the process for quasi-static loading steps.

\subsection{Combining Anderson acceleration and over-relaxation}
\begin{figure}
\centering
\includegraphics[width = 0.5\textwidth]{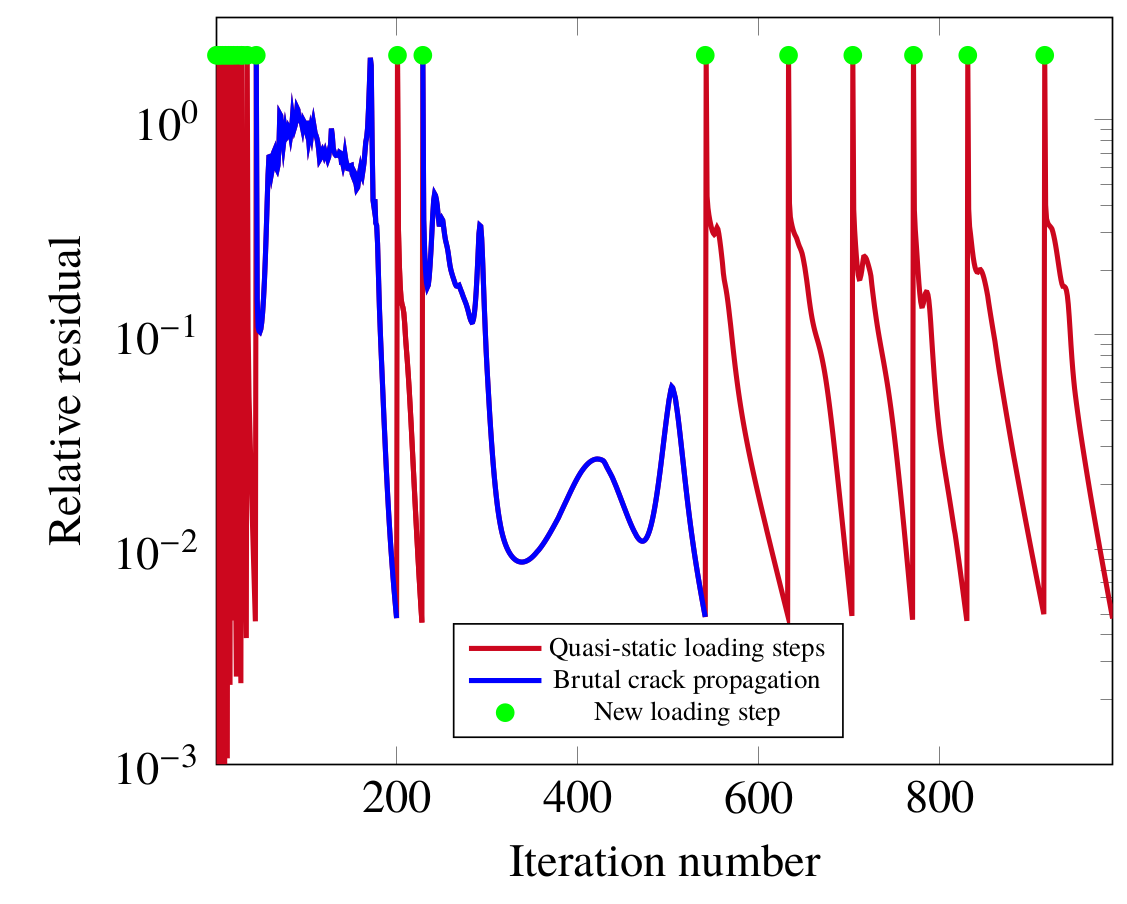}
\caption{{\bf Asymmetrical bending test:} Relative residual evolution (see equation~\eqref{eq:relresidual}) over the simulation. See Section~\ref{sec:holes} for an explanation of the test case. Similar behavior is experienced for all proposed test cases.}
\label{fig:motivation}
\end{figure}

As neither Anderson acceleration nor over-relaxation should be applied naively to the staggered scheme \eqref{eq:solutionalternate1}--\eqref{eq:solutionalternate2}, due to their mentioned weaknesses, we propose a combined robust acceleration scheme. The key observations that motivate such a method are:
\begin{itemize}
 \item Anderson acceleration is locally accelerating, while over-relaxation might struggle close to the solution.
 \item Anderson acceleration is applied as a post-processing algorithm to the increments of the staggered scheme, hence, only trivial implementation is required to switch between relaxation and Anderson acceleration. 
 \item During crack propagation the residuals for the staggered scheme show a stagnating, oscillatory behavior, and during quasi-static steps they are strictly decreasing, see \cite{gerasimovline} and Figure~\ref{fig:motivation}. Therefore, it is possible to use residual evolution as a rule for switching between the acceleration techniques. 
\end{itemize}

 A new parameter $N_{\omega\rightarrow AA}\in \mathbb{N}$, related to the switch from relaxation to Anderson acceleration is defined, and at loading step $n$ the combined accelerated staggered scheme reads: 
\begin{enumerate}
 \item Apply with Anderson acceleration of given depth $m$.
 \item While the norms of the residuals are strictly decreasing, continue with Anderson acceleration until convergence.
 \item If the norms of the residuals are not strictly decreasing, switch to relaxation with given parameter $\omega$.
 \item When the norms of the $N_{\omega\rightarrow AA}$ previous residuals are strictly decreasing go back to $1$, and restart\footnote{\label{fn:restart}Restart means that Anderson acceleration should be applied as it is in the first iteration, i.e., using no information of previous increments and iterates.} Anderson acceleration.
\end{enumerate}

Below, we give a pseudo-code for the new combined acceleration method. Define the residual norm $\mathrm{Res}_i = \left\|\mathrm{Res}_{\bm u}\left(\bm u_h^{n,i},\varphi_h^{n,i}\right)\right\|_2$ as in \eqref{eq:stopresabs}, and notice that the application of Anderson acceleration and relaxation in the pseudo-code denotes the $i$-th step of the accelerations (see Algorithm~\ref{alg:AA} and Algorithm~\ref{alg:OR}).\\
\begin{algorithm}[H]
Given depth $m$, relaxation $\omega$, initial guess $\varphi^{n,0}_h$, and switch $N_{\omega\rightarrow AA}$\;
\it{relaxing} := \it{False}\; 
\For{$i = 1,2,...$ until convergence}{
\If{not(relaxing)}{
\If{$i=1$ or $\mathrm{Res}_{i}\leq \mathrm{Res}_{i-1}$}{
apply Anderson accelerated staggered scheme, giving $(\bm u_h^{n,i}, \varphi^{n,i}_h)$\;}
\Else{
apply relaxed staggered scheme, giving $(\bm u_h^{n,i}, \varphi^{n,i}_h)$\;
\it{relaxing} := \it{True}\; 
}
}
\Else{
\If{not($\mathrm{Res}_{i}\leq \mathrm{Res}_{i-1}\leq\dots\leq \mathrm{Res}_{i-N_{\omega\rightarrow AA}-1}$)}{
apply relaxed staggered scheme, giving $(\bm u_h^{n,i}, \varphi^{n,i}_h)$\;
}
\Else{
$\mathrm{\it{restart}}^1$ and apply Anderson accelerated staggered scheme, giving $(\bm u_h^{n,i}, \varphi^{n,i}_h)$\;
\it{relaxing} := \it{False}\; 
}
}
}
\caption{Combined algorithm}
\label{alg:comb}
\end{algorithm}

\section{Numerical examples}\label{sec:numerics}
This section explores the effects of the proposed acceleration methods from Section~\ref{sec:acceleration} applied to the staggered scheme \eqref{eq:solutionalternate1}--\eqref{eq:solutionalternate2}. Both Anderson acceleration and over-relaxation alone are shown to be infeasible acceleration methods when plainly applied to the staggered scheme, while the combined scheme is superior to the unaccelerated scheme for all tests.
We consider four different test cases which are widely used for numerical studies in the literature:
\begin{itemize}
 \item A domain with a single notch subject to
 \begin{itemize}
  \item tensile load;
  \item shear load.
 \end{itemize}
 \item An L-shaped domain subject to loading.
 \item Bending of an asymmetrically notched beam with holes.

\end{itemize}
All the numerical examples have been implemented using modules from the DUNE project \cite{dune2015}, specifically dune-functions \cite{dunefunctions2015,dunefunctions2018}.

\subsection{Single notch test}
\begin{figure}
\begin{subfigure}{0.45\linewidth}
\includegraphics[width=\textwidth]{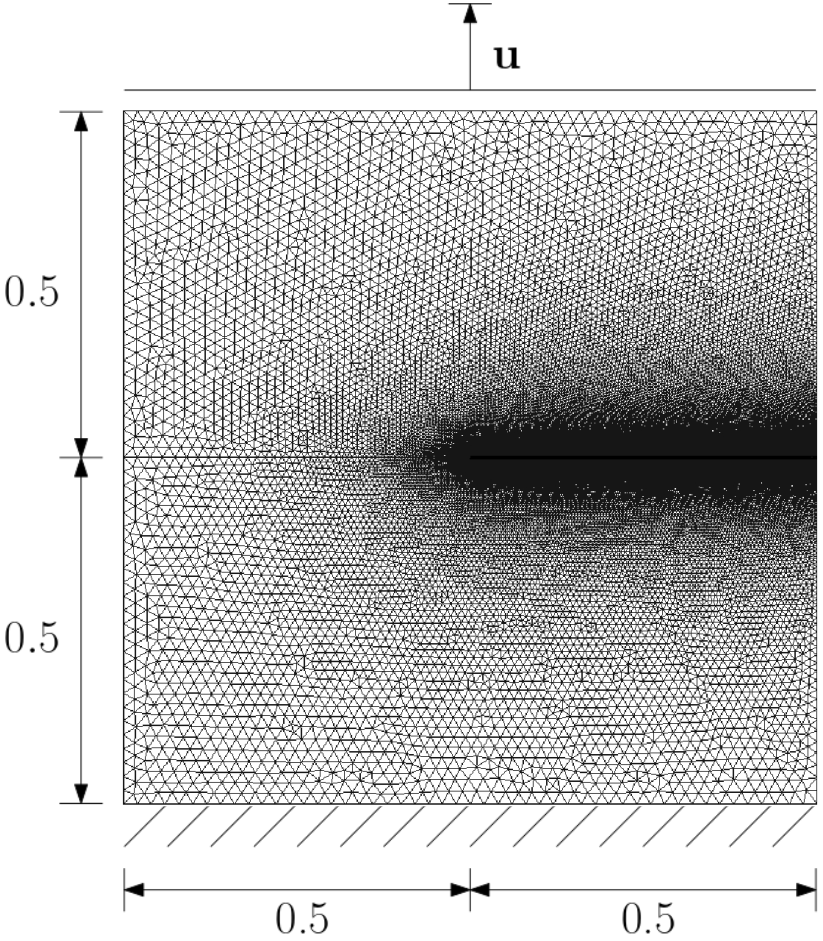}
\caption{{\bf Single notch tensile test:} The bottom boundary of the domain is fixed ($\bm u={\bf 0}$), and the top boundary is uniformly displaced over time in the vertical direction ($u_y = \bar{u}n$) while fixed in the horizontal direction ($u_x=0$). The mesh is refined according to the expected crack path and contains a total of 36995 nodes.\vspace{\baselineskip}}
\label{fig:notch-tensile-domain}
\end{subfigure}
\hspace{0.10\linewidth}
\begin{subfigure}{0.45\linewidth}
\includegraphics[width=\textwidth]{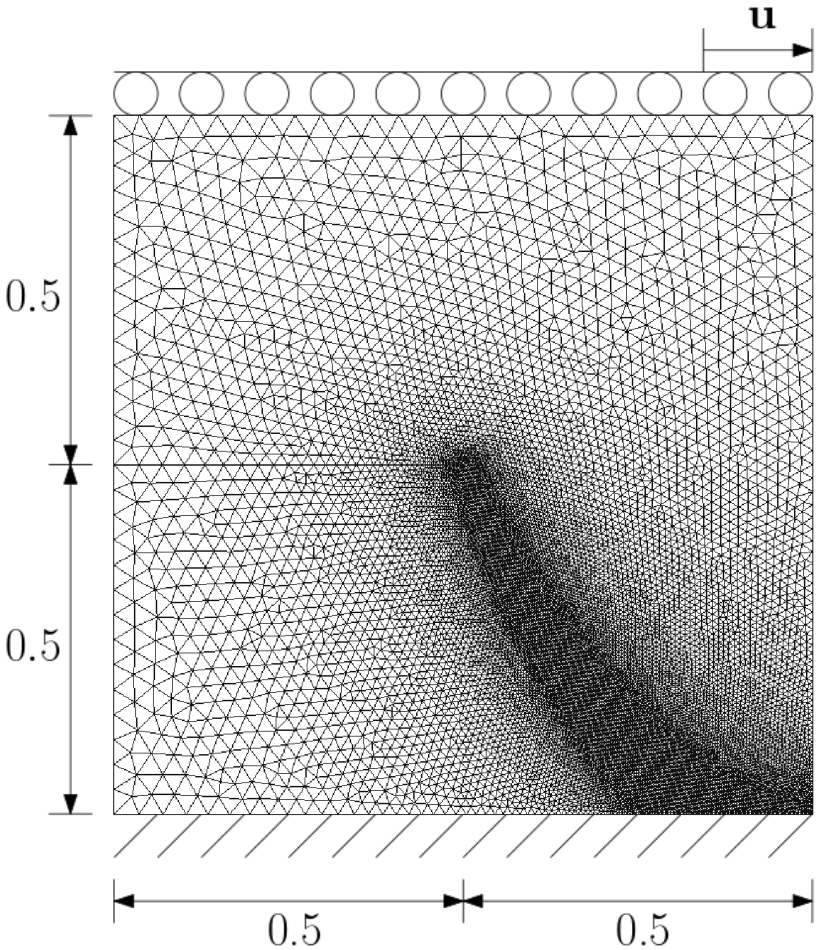}
\caption{{\bf Single notch shear test:} The bottom boundary is fixed ($\bm u={\bf 0}$), and the top boundary is uniformly displaced over time in the horizontal direction ($u_x = \bar{u}n$). The left, right and top boundaries, and the lower lip of the prescribed crack are fixed in the vertical direction ($u_y=0$). The mesh is refined according to the expected crack path and contains a total of 12660 nodes.}
\label{fig:notch-shear-domain}
\end{subfigure}
\caption{Domain, boundary conditions, and mesh for the single-edge notch test cases.}
\label{fig:notch-domain}
\end{figure}

Two of the most commonly found test cases in the literature are both based on the same single notch geometry \cite{miehethermodynamics, ambati}. They consist of a square domain with a pre-existing crack that penetrates half the domain, see Figure~\ref{fig:notch-domain}. The domain is held still at the bottom, and a displacement driven load is applied at the top boundary.

\subsubsection{Single notch tensile test case}
A tensile load is applied on the top boundary, and at loading step $n$ we have 
$$\bm u^n_{|_{\Gamma^{\mathrm{Top}}}}=\begin{pmatrix}0\\ \bar{u}n \end{pmatrix},$$ where the load size $\bar{u}$ is given in Table~\ref{tab:notch-test} and $\Gamma^{\mathrm{Top}}$ is the top part of the boundary in Figure~\ref{fig:notch-tensile-domain}. Due to the load being strictly tensile, there is no need to split the elastic strain energy functional \eqref{eq:elasticenergy} into tensile and compressive parts, which would effectively add nonlinearities to the system. Therefore, the first term in \eqref{eq:variationalmechanics} is replaced by $\left\langle g(\varphi)\bm \sigma(\bm u)\!:\! \bm\varepsilon(\bm u))\right\rangle,$ 
for \begin{equation}\label{eq:load}\bm \sigma(\bm u ) := \dfrac{\partial \bm \psi(\bm\varepsilon(\bm u))}{\partial\bm\varepsilon(\bm u)}.\end{equation} Material parameter values are chosen as in e.g., \cite{miehethermodynamics}, and can be found in Table \ref{tab:notch-test}. We employ a triangular mesh, which has been locally refined in the region where the crack is expected to propagate, see Figure~\ref{fig:notch-tensile-domain}.

\begin{figure}
\begin{subfigure}{0.45\linewidth}
\includegraphics[width=\textwidth]{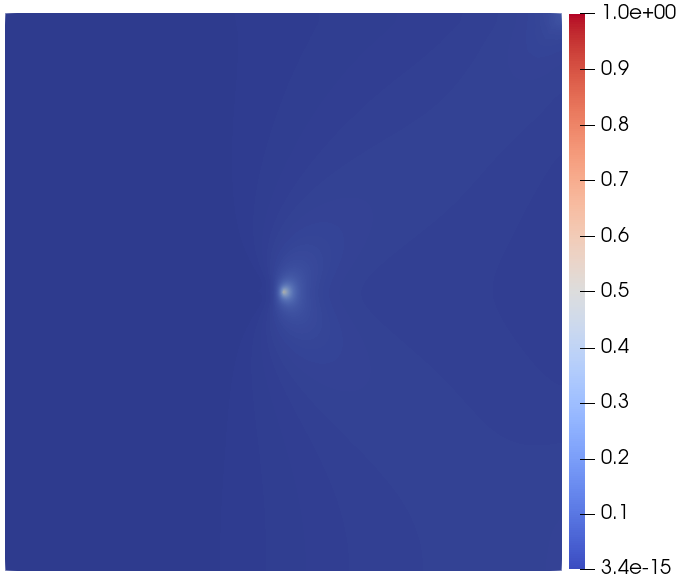}
\caption{Solution before crack growth at loading step 27.}
\label{fig:sol-tensile-27}
\end{subfigure}
\hspace{0.10\linewidth}
\begin{subfigure}{0.45\linewidth}
\includegraphics[width=\textwidth]{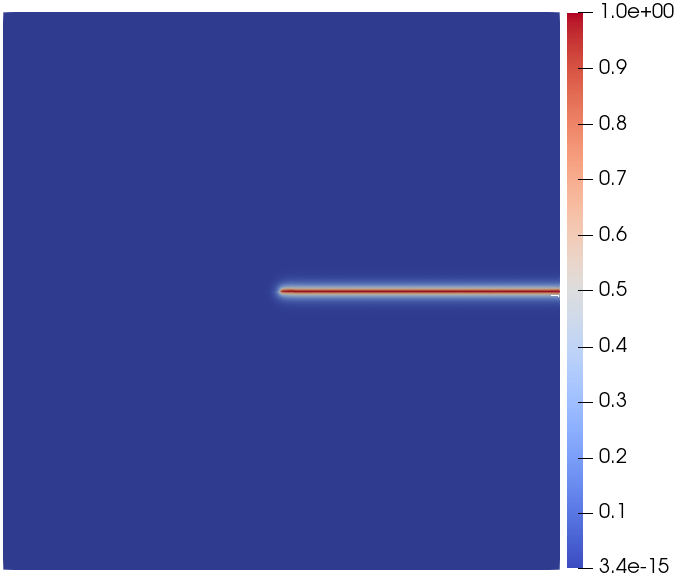}
\caption{Solution after crack propagation at loading step 28. }
\label{fig:sol-tensile-28}
\end{subfigure}
\caption{Solution for $\varphi$ for the single notch tensile test case.}
\label{fig:tensile-solution}
\end{figure}

\begin{figure}
\begin{subfigure}{0.45\linewidth}
 \includegraphics[width = \linewidth]{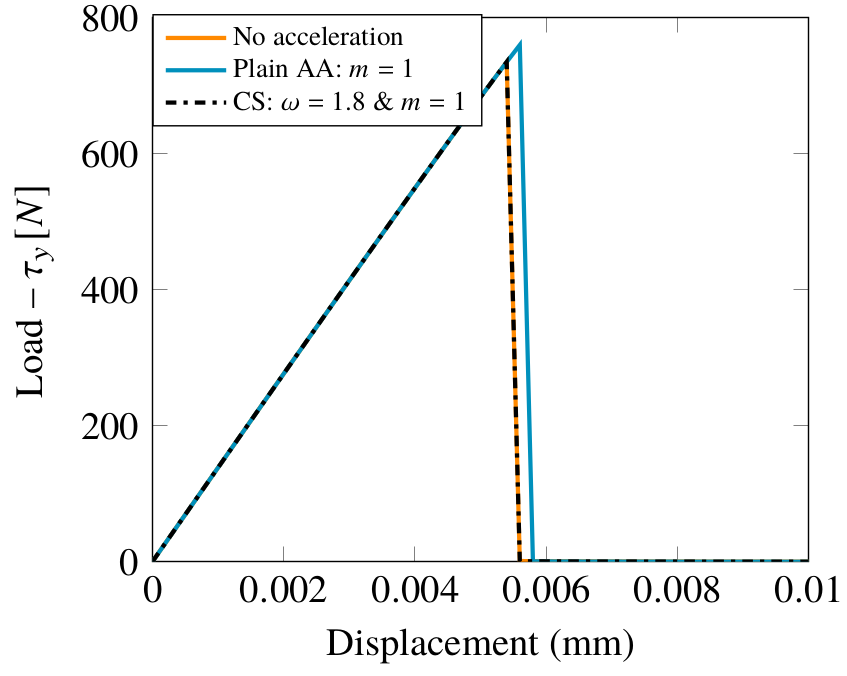}
\caption{Load curves.}
\label{fig:tensile-load}
\end{subfigure}
\hspace{0.10\linewidth}
\begin{subfigure}{0.45\linewidth}
\includegraphics[width=\linewidth]{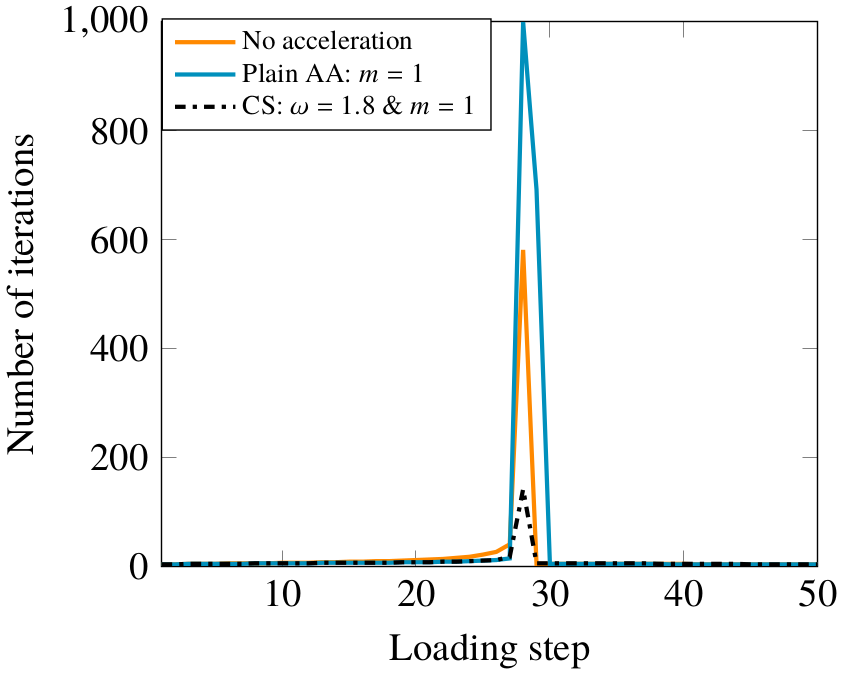}
\caption{Number of iterations per loading step.}
\label{fig:tensile-iterations}
\end{subfigure}
\caption{{\bf Single notch tensile test:} Load curves and number of iterations per loading step. ``AA'' is an abbreviation of Anderson acceleration, and ``CS``\ is the combined acceleration scheme. ''Plain AA'' means that Anderson acceleration is applied without any form of safeguard or combination with relaxation.}
\end{figure}

In this test case, the crack fully propagates in one single critical loading step, see Figure~\ref{fig:tensile-solution}, in which the crack gradually expands through the domain with increasing staggered iteration count. Figure~\ref{fig:tensile-iterations} shows that, as expected, the staggered scheme under Anderson acceleration alone struggles as a consequence of its local convergence. The combined scheme, however, takes advantage of over-relaxation and its ability to move further each iteration and accelerates this particular loading step significantly. For the rest of the loading steps, the combined scheme accelerates by Anderson acceleration, as its local convergence is sufficient. The total number of iterations for the combined scheme is, therefore, smaller than those of the unaccelerated staggered scheme and the Anderson accelerated staggered scheme.  

The traction vector is defined by 
\begin{equation}\label{eq:traction}\bm \tau = (\tau_x,\tau_y) = \int_{\Gamma^\mathrm{Top}}\bm\sigma\cdot \bm \nu \;dS,
\end{equation} 
where $\bm \nu$ is the outward pointing normal vector, and $\bm \sigma$ is defined in \eqref{eq:load}. For this problem, the load in the direction of interest is $\tau_y$, and we observe in Figure~\ref{fig:tensile-load} that the load-displacement curves remain unchanged after the combined acceleration. This is an important observation that demonstrates that the acceleration method only affects the convergence properties of the solver, not the quality of the solution. The Anderson accelerated staggered scheme, however, does not converge in the maximal prescribed iterations for each loading step and we observe that its load-displacement curve is affected. 

\begin{figure}
\includegraphics[width = \textwidth]{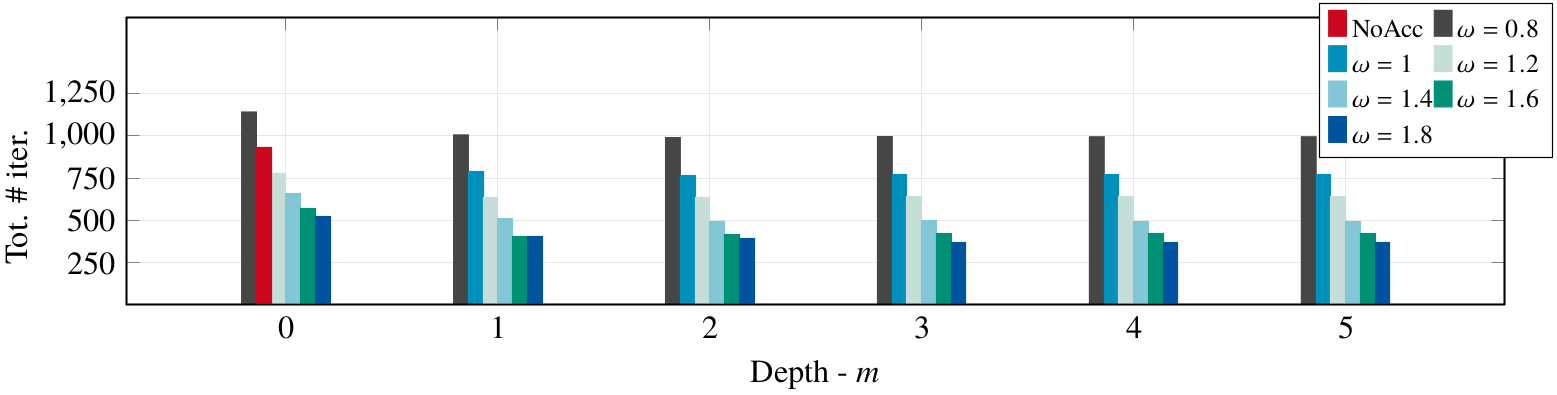}
 \caption{{\bf Single notch tensile test:} Total number of iterations for the combined scheme with different relaxation parameters and Anderson acceleration depths. ``NoAcc'' is the unaccelerated staggered scheme.}
 \label{fig:tensile-totaliterations}
\end{figure}

The total number of iterations is displayed in Figure~\ref{fig:tensile-totaliterations}, and there are several key observations. First of all, the combined scheme accelerates by more than 50 \% for large relaxation parameters. We also observe that the depth of Anderson acceleration is not influential as long as it is larger than one. Moreover, there is a trend that more aggressive over-relaxation (higher $\omega$) results in faster computations. Additionally, the combined acceleration scheme is robust with respect to the tuning parameters, $m$ and $\omega$, and exhibits convergence for all tested combinations.

\begin{table}[ht]
\begin{center}
\begin{tabular}{l|l|l|l}
 {\bf Parameter} & {\bf Symbol} &{\bf Value -- Tensile} &{\bf Value -- Shear}\\
 \hline 
 Lam\'e's 1. parameter&$\lambda$ & $121.15$ kN/$\mathrm{mm}^2$ & $121.15$ kN/$\mathrm{mm}^2$ \\
 Lam\'e's 2. parameter&$\mu $ & 80.77 kN/$\mathrm{mm}^2$ & 80.77 kN/$\mathrm{mm}^2$ \\
 Regularization width &$\ell$ & 0.0075 mm & 0.0075 mm \\
 Griffith's constant&$\mathcal{G}_c$ & $2.7$ N/mm & $2.7$ N/mm \\
 Tot. \# loading steps & N & 50 & 150 \\
 Load size &$\bar{u}$ & $2\cdot10^{-4}$ mm&$10^{-4}$ mm\\
 Fine mesh size & $h$ & 0.001 mm & 0.00375 mm\\
 Min. relax. steps & $N_{\omega\rightarrow \mathrm{AA}}$ & 5 & 5\\
 Abs.\ tol. & $\mathrm{Tol}_{\mathrm{Res/Inc,Abs}}$ & $10^{-8}$ & $10^{-8}$\\
 Rel.\ residual tol. & $\mathrm{Tol}_{\mathrm{Res,Rel}}$ &$5\cdot10^{-3}$ &$5\cdot10^{-3}$\\
Rel.\ increment tol. & $\mathrm{Tol}_{\mathrm{Inc,Rel}}$ &$10^{-2}$ &$10^{-2}$\\
Max.\ iter.\ pr.\ load.\ step & $\mathrm{Max}_\mathrm{iter}$ & 1000 & 1000\\
 Inner Newton tol. & $\mathrm{Tol}_{\mathrm{Inner}}$ &$10^{-4}$ &$10^{-4}$\\
\end{tabular}
\caption{Parameter values for the single notch test cases.}
\label{tab:notch-test}
\end{center}
\end{table}

\begin{remark}\label{rem:plain}
Notice that the plots of the number of iterations for depth $m= 0$ in Figures~\ref{fig:tensile-totaliterations}, \ref{fig:shear-totaliterations}, \ref{fig:lshape-totaliterations} and \ref{fig:holes-totaliterations} are not corresponding to plain relaxation. Here, relaxation is switched on and off depending on residual evolution, turning it into a safeguarded relaxation. The same goes for the plots of plain Anderson acceleration with over-relaxation parameter $\omega=1$. These correspond to safeguarded Anderson accelerations, similar to those that are proposed in \cite{AAsafetyguards}.
\end{remark}

\subsubsection{Single notch shear case}
In this test case, a shear load is applied on the top boundary of a unit square domain with a prescribed crack that halfway penetrates the domain. The displacement boundary condition
$$\bm u^n_{|_{\Gamma^{\mathrm{Top}}}}=\begin{pmatrix} \bar{u}n \\0\end{pmatrix}$$ is applied at loading step $n$. The load size $\bar{u}$ is presented in Table~\ref{tab:notch-test}, and the top part of the boundary $\Gamma^{\mathrm{Top}}$ is displayed together with more details on the domain and boundary conditions in Figure~\ref{fig:notch-shear-domain}.
The material properties are taken from \cite{miehethermodynamics} and displayed in Table~\ref{tab:notch-test}. A triangular mesh, which is refined according to where the crack is expected to propagate, has been employed, see Figure~\ref{fig:notch-shear-domain}. 

\begin{figure}
\begin{subfigure}{0.22\linewidth}
\includegraphics[width=\textwidth]{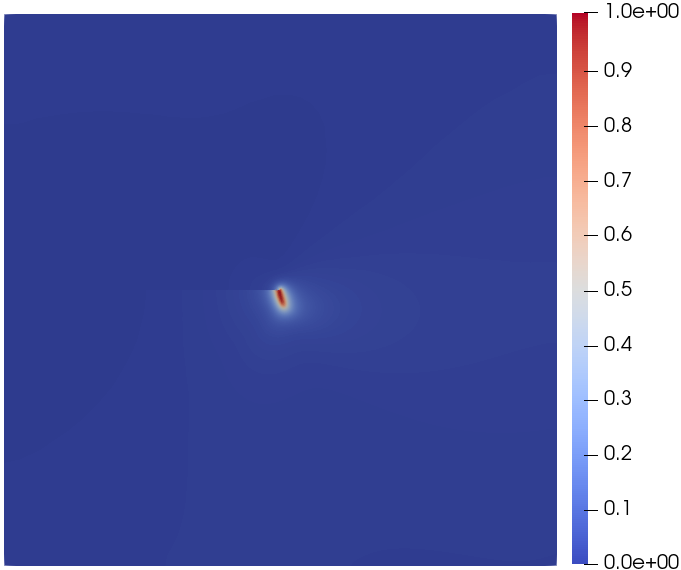}
\caption{Solution at loading step 100.}
\label{fig:sol-shear-100}
\end{subfigure}
\hspace{0.03\linewidth}
\begin{subfigure}{0.22\linewidth}
\includegraphics[width=\textwidth]{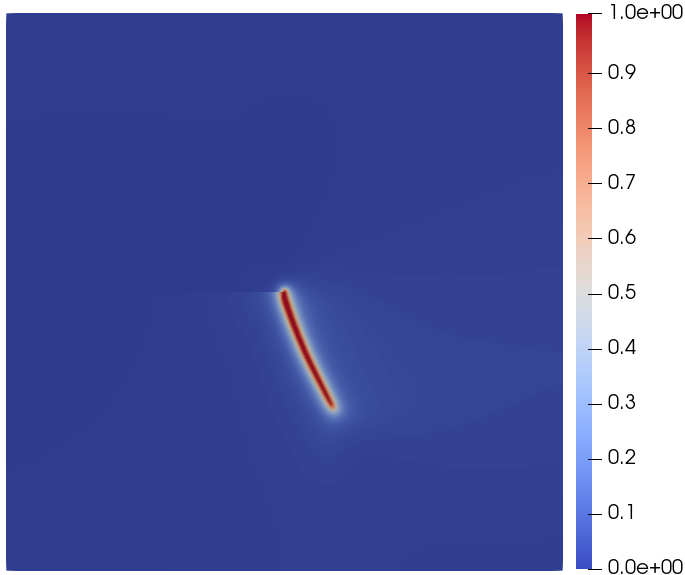}
\caption{Solution at loading step 117. }
\label{fig:sol-shear-150}
\end{subfigure}
\hspace{0.03\linewidth}
\begin{subfigure}{0.22\linewidth}
\includegraphics[width=\textwidth]{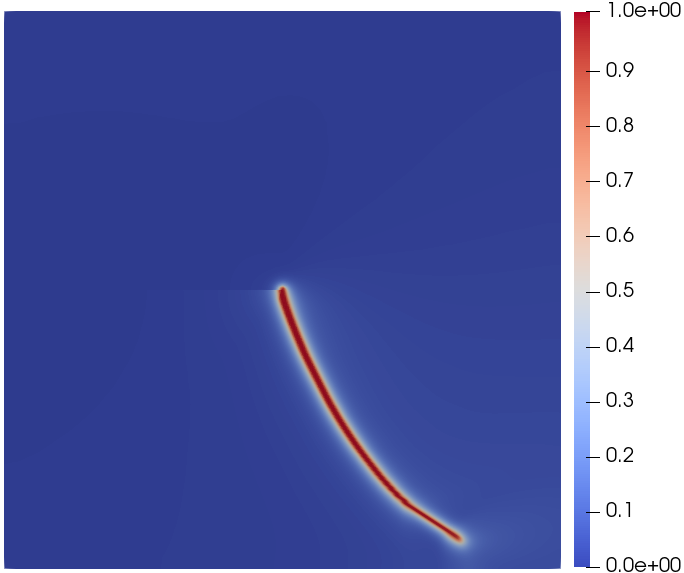}
\caption{Solution at loading step 133. }
\label{fig:sol-shear-150}
\end{subfigure}
\hspace{0.03\linewidth}
\begin{subfigure}{0.22\linewidth}
\includegraphics[width=\textwidth]{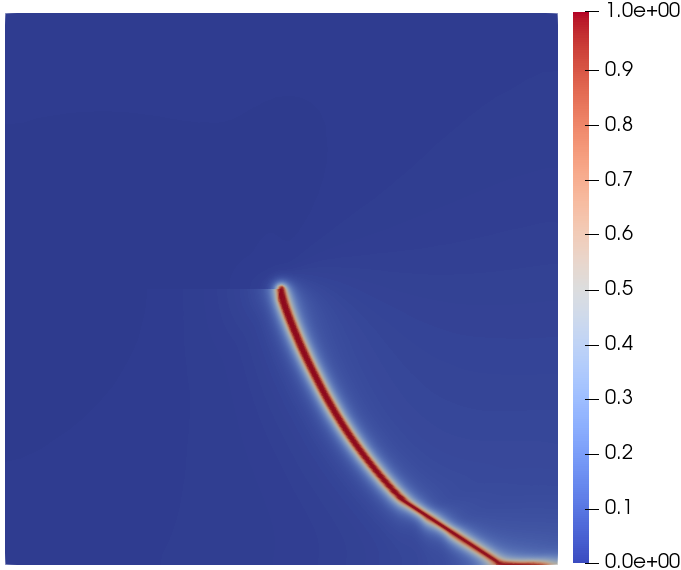}
\caption{Solution at loading step 150. }
\label{fig:sol-shear-150}
\end{subfigure}
\caption{Solution for $\varphi$ for the single notch shear test case.}
\label{fig:sol-shear}
\end{figure}

\begin{figure}
\begin{subfigure}{0.45\linewidth}
 \includegraphics[width = \textwidth]{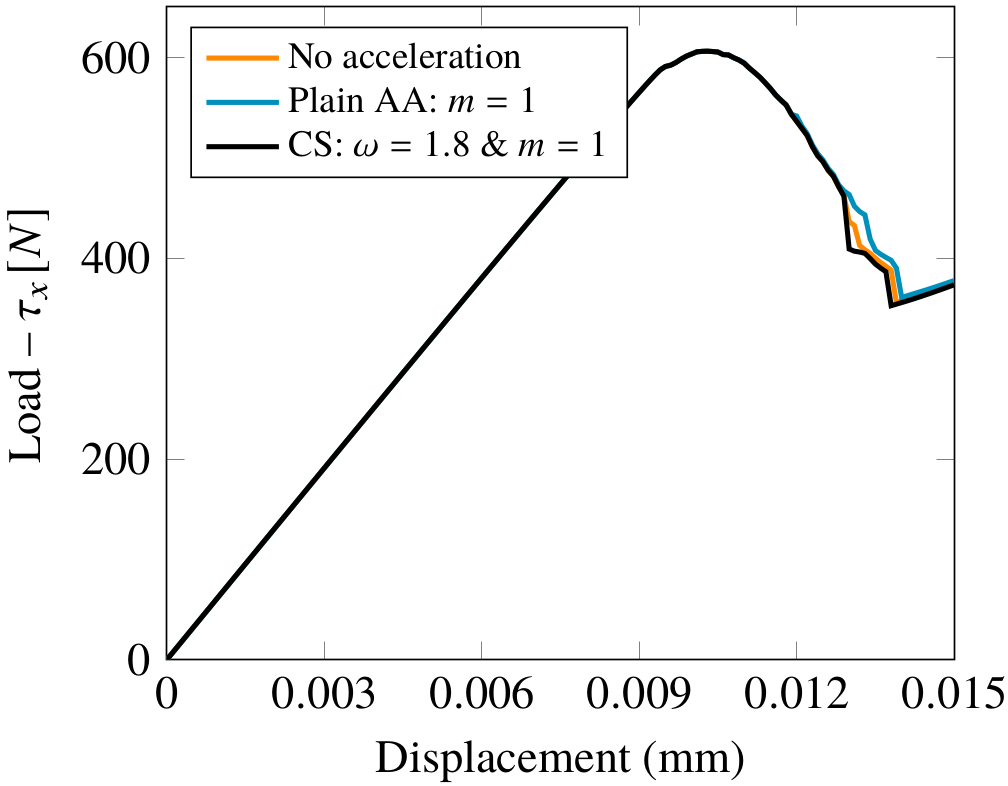}
\caption{Load curves.}
\label{fig:shear-load}
\end{subfigure}
\hspace{0.1\linewidth}
\begin{subfigure}{0.45\linewidth}
 \includegraphics[width = \textwidth]{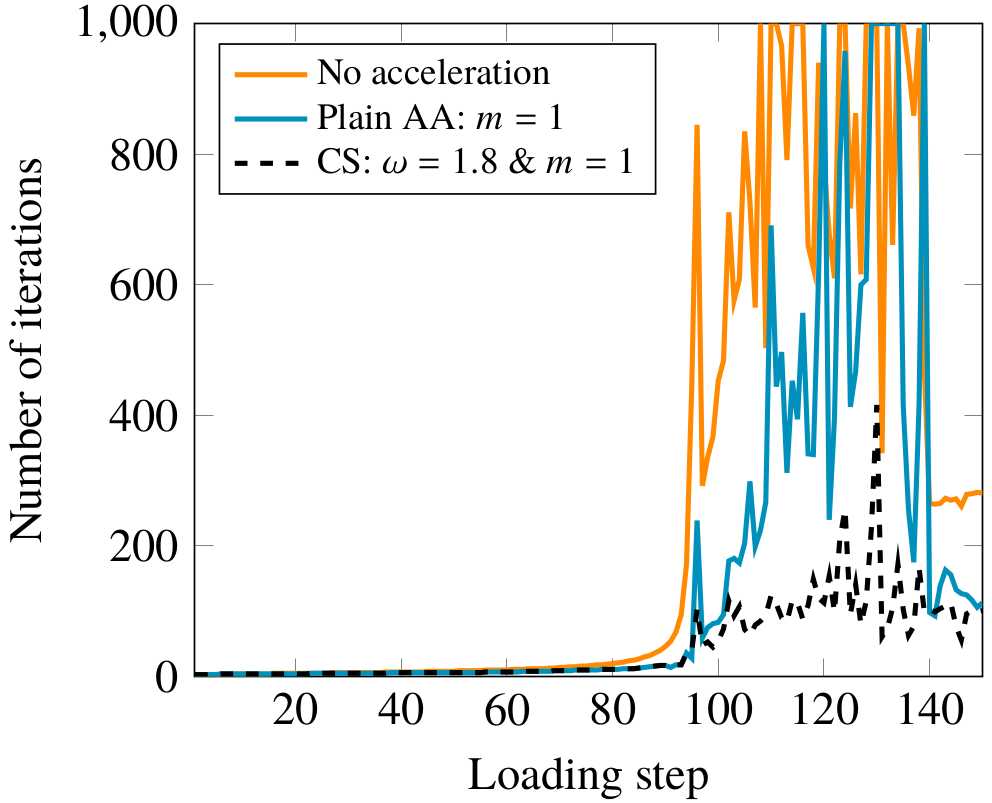}
\caption{Number of iterations per loading step. Maximum number of iterations per loading step is 1000.}
\label{fig:shear-iterations}
\end{subfigure}
\caption{{\bf Single notch shear test:} Load curves and number of iterations per loading step. ``AA'' is an abbreviation of Anderson acceleration, and ``CS``\ is the combined acceleration scheme. ''Plain AA'' means that Anderson acceleration is applied without any form of safeguard or combination with relaxation.}
\end{figure}

Contrary to the tensile test case, the crack propagation happens gradually over the course of many loading steps, see Figure~\ref{fig:sol-shear}. Therefore, solutions at subsequent loading steps do not differ as significantly as for the brutal crack growth in the tensile test case. We expect that the Anderson acceleration is a more suitable choice for accelerating the staggered scheme. Indeed, Figure~\ref{fig:shear-iterations} shows that even with the naive Anderson acceleration the staggered scheme is quite significantly accelerated. Moreover, the combined scheme is even better, and we see that it reaches convergence in every single loading step. 

In the load-displacement curves, Figure~\ref{fig:shear-load}, the load $\tau_x$ from \eqref{eq:traction} is displayed for each loading step. The plot shows minor differences towards the end of the displacement. This is due to the scheme not converging in its given maximal iterations per loading step (see Table~\ref{tab:notch-test}) for both the unaccelerated staggered scheme and the plain Anderson accelerated scheme in all the loading steps. This is similar to the tensile case where the load-displacement curves (Figure~\ref{fig:tensile-load}) also are affected.

\begin{figure}
  \includegraphics[width = \textwidth]{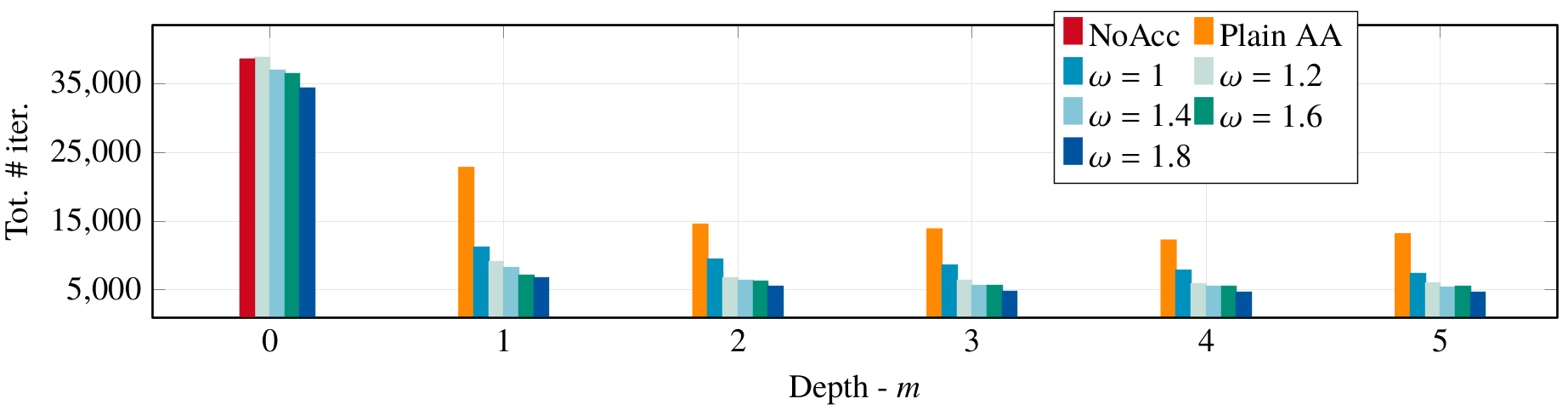}
 \caption{{\bf Single notch shear test:} Total number of iterations for different relaxation parameters and Anderson acceleration depths. ``NoAcc'' is the unaccelerated staggered scheme, and ``Plain AA`` is Anderson acceleration without the combination with relaxation.}
 \label{fig:shear-totaliterations}
\end{figure}

In Figure~\ref{fig:shear-totaliterations}, we see the total iteration count for several acceleration depths in combination with over-relaxation. The figure shows that the staggered scheme is accelerated significantly for all combinations of Anderson acceleration and over-relaxation as long as the depth is greater than one. In fact, we have more than 80 \% reduction in the total number of iterations when choosing a high relaxation parameter. Moreover, for this test case the plain Anderson acceleration is in itself a suitable alternative to the unaccelerated staggered scheme. Notice the difference between plain Anderson acceleration and Anderson acceleration combined with over-relaxation of depth one described in Remark~\ref{rem:plain}.

\subsection{L-shaped domain subject to loading}
\begin{figure}
\begin{subfigure}{0.45\linewidth}
\includegraphics[width=0.9\textwidth]{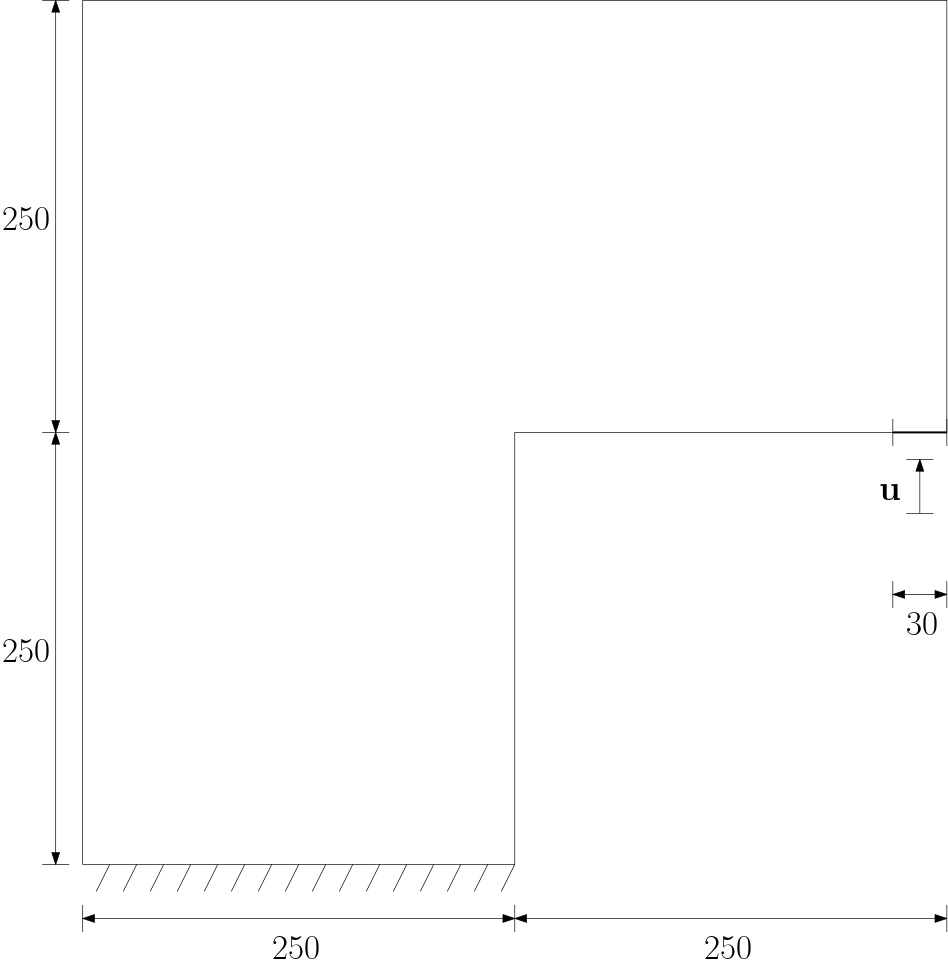}
\caption{{\bf L-shape test:} Domain and boundary conditions. The bottom boundary of the domain is fixed ($\bm u={\bf 0}$), and the boundary-segment $\Gamma_D = \{(x,y)\in \Omega\;|\; 470\leq x\leq500, \; y = 250\}$ ($(0,0)$ is in the lower left corner) is uniformly displaced over time in the vertical direction ($u_y = \bar{u}n$).\vspace{2em}}
\label{fig:lshape-domain}
\end{subfigure}
\hspace{0.10\linewidth}
\begin{subfigure}{0.45\linewidth}
\includegraphics[width=\textwidth]{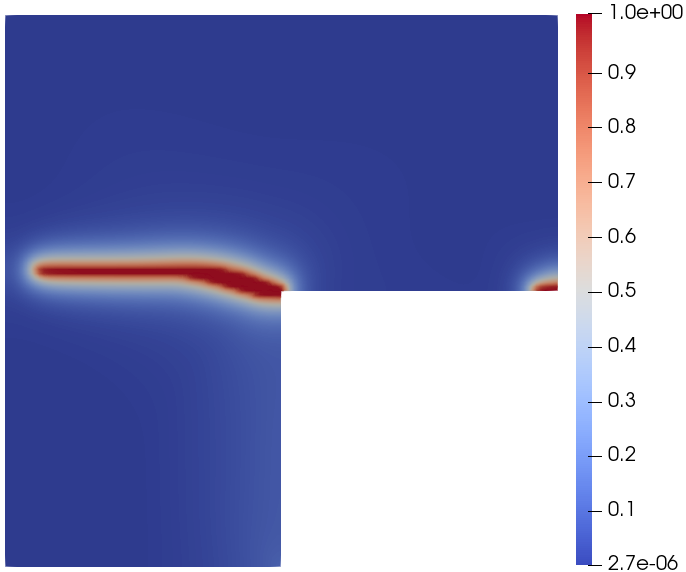}
\vspace{1em}
\caption{{\bf L-shape test:} The solution for $\varphi$ for the L-shape test after 800 loading steps.\vspace{6.4em}}
\label{fig:lshape-solution}
\end{subfigure}
\caption{Domain with boundary conditions and solution for the L-shape test.}
\end{figure}

An L-shaped domain with a displacement boundary condition applied on the right part of the boundary is considered, see Figure~\ref{fig:lshape-domain} for details. The displacement is uniformly increased on the boundary segment over 800 loading steps. As a result, a crack occurs in the inner corner, propagating into the domain, see Figure~\ref{fig:lshape-solution}. A uniform quadrilateral mesh with a mesh diameter of $\frac{125}{32}$ mm is employed. See Table~\ref{tab:holes} for material and computational parameters.

\begin{figure}
\begin{subfigure}{0.45\linewidth}
\includegraphics[width=\textwidth]{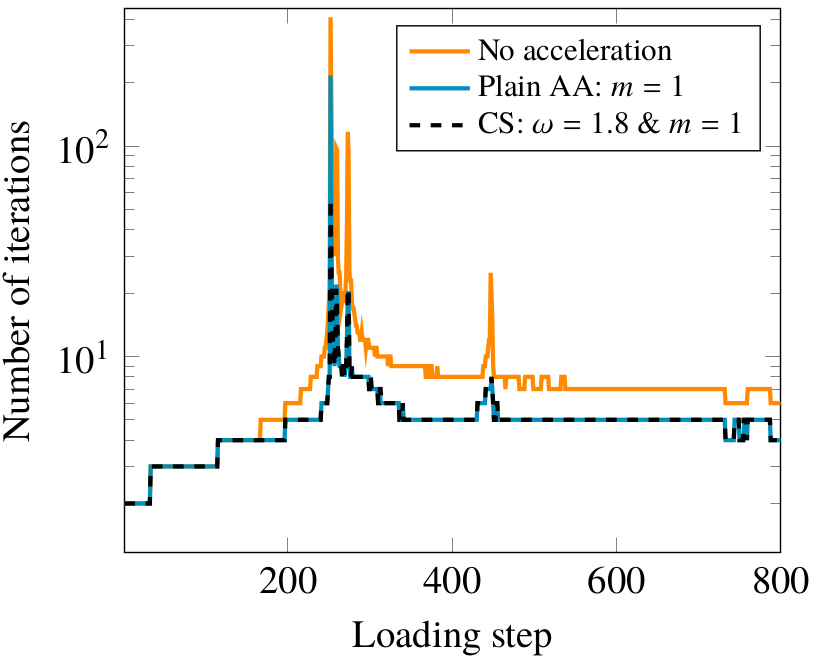}
\caption{Number of iterations per loading step. `AA'' is an abbreviation of Anderson acceleration, and ``CS``\ is the combined acceleration scheme. ''Plain AA'' means that Anderson acceleration is applied without any form of safeguard or combination with relaxation. Notice that the plot has a log-scale on the y-axis.}
\label{fig:lshape-iterations}
\end{subfigure}
\hspace{0.1\linewidth}
\begin{subfigure}{0.45\linewidth}
\includegraphics[width=\textwidth]{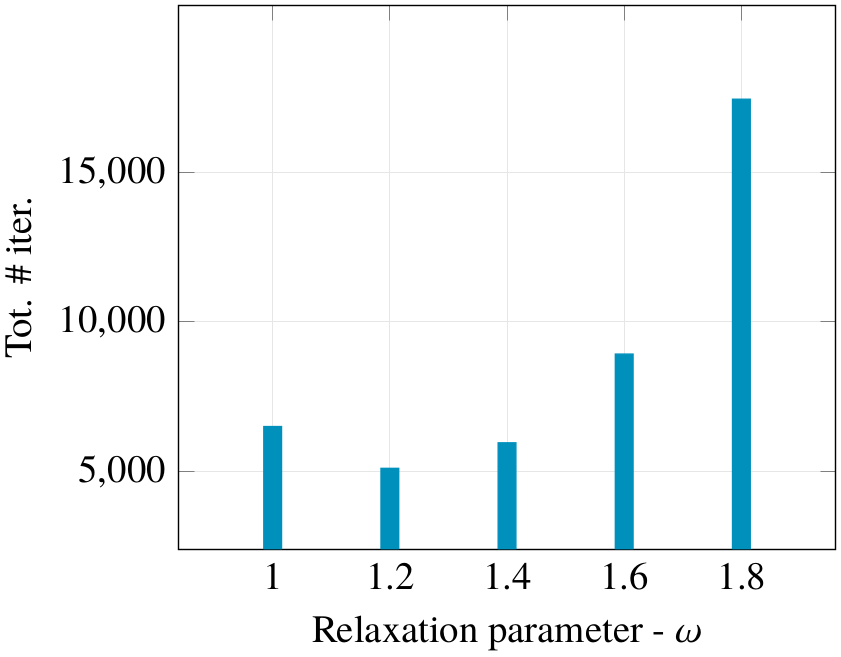}
 \caption{Total number of iterations for different over relaxations parameters applied without safeguard or combination.\vspace{4\baselineskip}}
 \label{fig:lshape-totaliterations-OR}
\end{subfigure}
\caption{{\bf L-shape test:} number of iterations per loading step and total iterations for several over-relaxation parameters.}
\end{figure}

Here, the crack propagation has a character somewhere between the single notch tensile test and the single notch shear test. Crack initiation shows similar behavior as brutal crack propagation, but not as extreme as for the single notch tensile test. A large peak in the number of iterations is experienced when the crack initiates, see Figure~\ref{fig:lshape-iterations}. We observe that both the combined scheme and plain Anderson acceleration with depth $m=1$ accelerates for all loading steps. Moreover, the only difference between the combined acceleration and Anderson acceleration is in the large peak where the combined scheme outperforms Anderson acceleration. For the rest of the simulation, the staggered scheme converges in relatively few iterations (less than 30 per loading step), but the accelerated method converges faster in almost every loading step. 

\begin{figure}
\includegraphics[width = \textwidth]{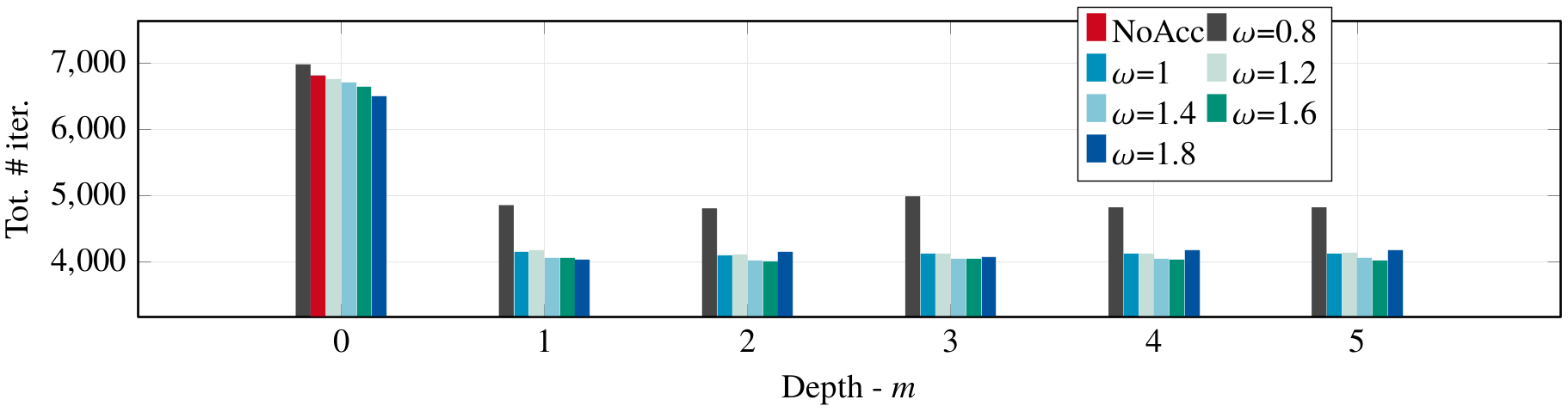}
\caption{{\bf L-shape test:} Total number of iterations for different relaxation parameters and Anderson acceleration depths. ``NoAcc'' is the unaccelerated staggered scheme.}
 \label{fig:lshape-totaliterations}
\end{figure}

Figure~\ref{fig:lshape-totaliterations-OR} displays the total number of iterations for plain over-relaxation with several relaxation parameters. A parabolic dependence on the relaxation parameter is observed, and choosing it to be too large results in more than three times the number of iterations that are required by the unaccelerated staggered scheme. Therefore, a plain application of over-relaxation is not recommended. The total number of iterations required by the combined acceleration, however, is significantly smaller than those of the unaccelerated staggered scheme, as observed in Figure~\ref{fig:lshape-totaliterations}. Although the reduction in the number of iterations is not as extreme as for the single notch shear test case the combined scheme accelerates robustly with respect to the tuning parameters. It is clear that any combination of Anderson acceleration and over-relaxation is superior to the unaccelerated staggered scheme, accelerating by approximately 40 \%.

\begin{figure}
 \begin{subfigure}{0.45\linewidth}
\includegraphics[width  =\textwidth]{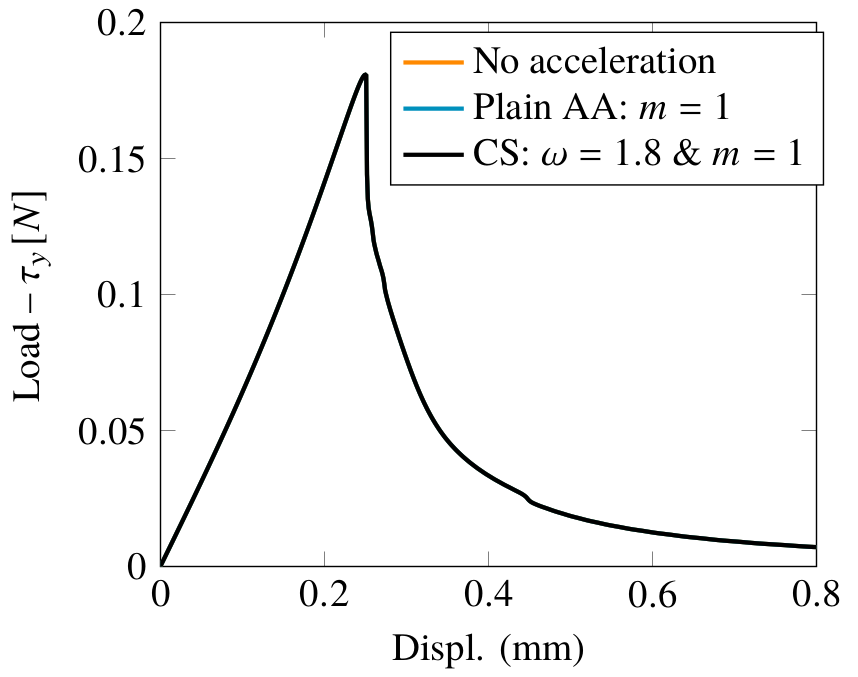}
\caption{{\bf L-shaped test} Load-displacement curve.}
\label{fig:lshape-load}
\end{subfigure}
\hspace{0.1\linewidth}
\begin{subfigure}{0.45\linewidth}
\includegraphics[width  =\textwidth]{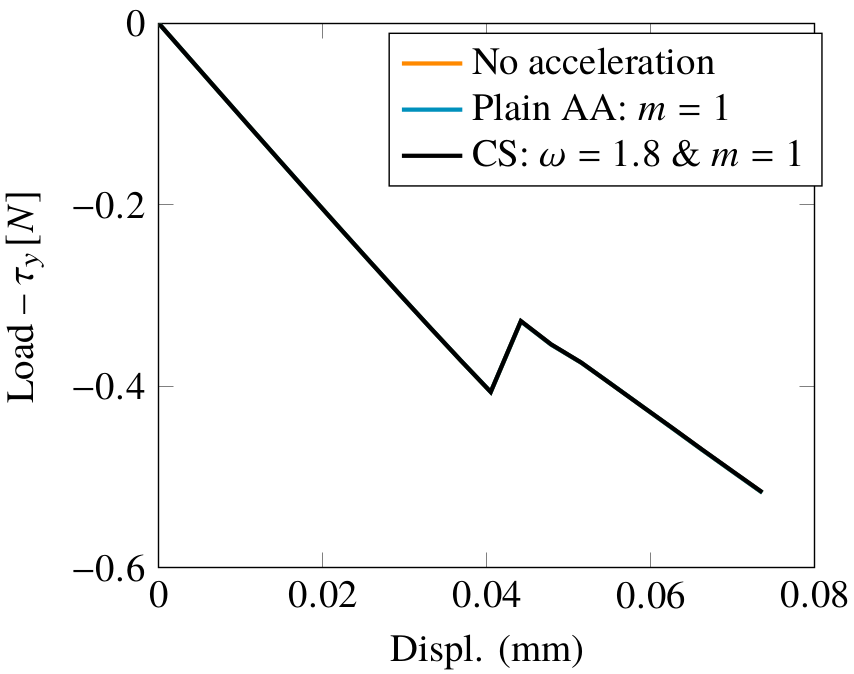}
\caption{{\bf Asymmetrical bending test:} Load-displacement curve.}
\label{fig:holes-load}
\end{subfigure}
\caption{``CS``\ is the combined acceleration scheme. ''Plain AA'' means that Anderson acceleration is applied without any form of safeguard or combination with relaxation.}
\end{figure}

The load-displacement curve for this test case is displayed in Figure~\ref{fig:lshape-load}. Here, the traction vector (see equation~\eqref{eq:traction}) is calculated on the bottom boundary and the vertical component $\tau_y$ is considered. We observe that, as all acceleration schemes converge within each loading step, the curves are completely overlapping.

\begin{table}[ht]
\begin{center}
\begin{tabular}{l|l|l|l}
 Parameter & Symbol & L-shape & Bend.\ test\\
 \hline 
 Lam\'e's 1. parameter&$\lambda$ & $6.16$ kN/$\mathrm{mm}^2$ & $8$ kN/$\mathrm{mm}^2$ \\
 Lam\'e's 2. parameter&$\mu $ & 10.95 kN/$\mathrm{mm}^2$ & 12 kN/$\mathrm{mm}^2$ \\
 Regularization width &$\ell$ & 10 mm & 0.1 mm \\
 Griffith's constant&$\mathcal{G}_c$ & $9.5\cdot10^{-5}$ kN/mm & $10^{-3}$ kN/mm \\
 Load size &$\bar{u}$ & $10^{-3}$ mm& $-10^{-2}e^{-\frac{(x-10)^2}{100}}$ mm\\
 Fine mesh size & $h$ & $\frac{125}{32}$ mm & 0.05 mm\\
Min. relax. steps & $N_{\omega\rightarrow \mathrm{AA}}$ & 5 & 5\\
 Abs. tol. & $\mathrm{Tol}_{\mathrm{Res/Inc,Abs}}$ & $10^{-8}$ & $10^{-8}$\\
 Rel. residual tol. & $\mathrm{Tol}_{\mathrm{Res,Rel}}$ &$5\cdot10^{-3}$ &$5\cdot10^{-3}$\\
Rel. increment tol. & $\mathrm{Tol}_{\mathrm{Inc,Rel}}$ &$10^{-2}$ &$10^{-2}$\\
Max.\ iter.\ pr.\ load.\ step & $\mathrm{Max}_\mathrm{iter}$ & 1000 & 1000\\
 Inner Newton tol. & $\mathrm{Tol}_{\mathrm{Inner}}$ &$10^{-4}$ &$10^{-4}$\\
\end{tabular}
\caption{Parameter values for the L-shape and asymmetrical bending tests.}
\label{tab:holes}
\end{center}
\end{table}

\subsection{Asymmetrical bending test }\label{sec:holes}
\begin{figure}
\begin{subfigure}{0.45\linewidth}
\includegraphics[width=\textwidth]{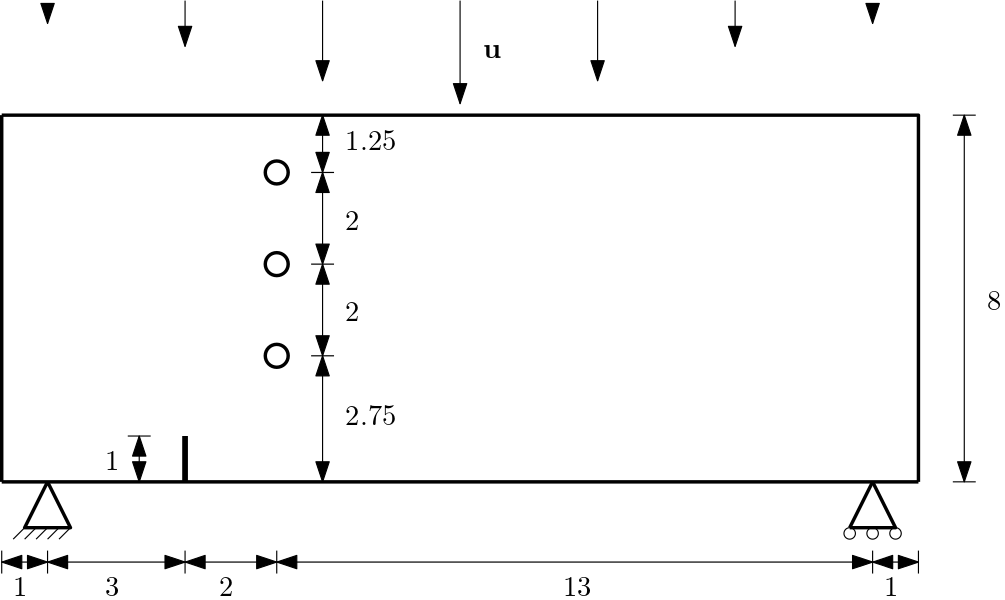}
\caption{{\bf Asymmetrical bending test:} Domain and boundary conditions. The three holes have 0.5 mm diameter. At the right and left boundaries, $\varphi = 0$ is enforced to prevent artificial crack initiation, see \cite{mesgarnejad}. The beam is fixed ($\bm u =\bm 0$) at the left foot and fixed in the vertical direction ($u_y=0$) at the right foot. Displacement condition $u_y=\bar{u}n$ is applied at the top boundary.}
\label{fig:holes-domain}
\end{subfigure}
\hspace{0.1\linewidth}
\begin{subfigure}{0.45\linewidth}
 \includegraphics[width=\textwidth]{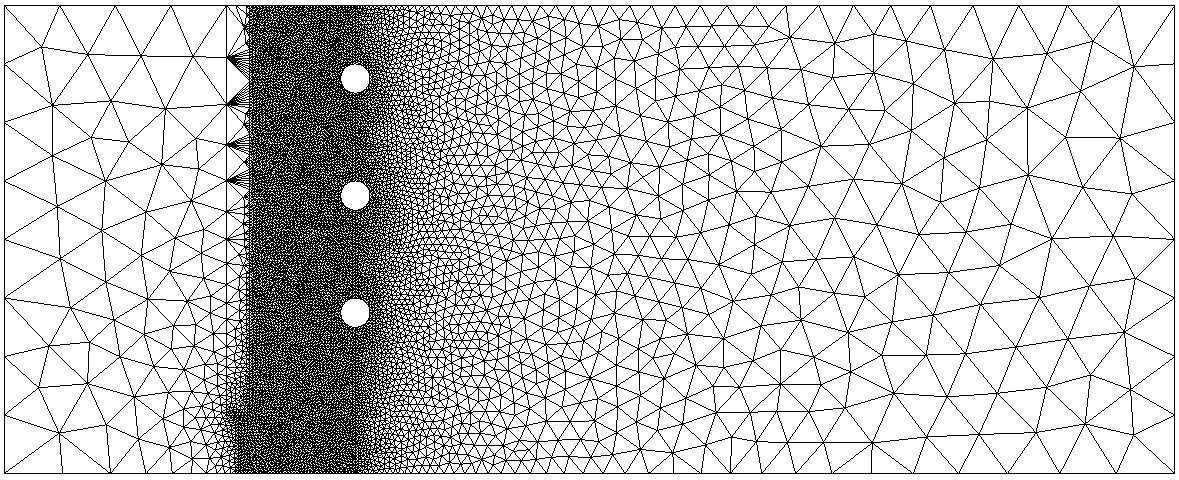}
 \vspace{0.5\baselineskip}
\caption{{\bf Asymmetrical bending test:} The mesh is refined according to the expected crack path and contains a total of 9598 nodes. 
\vspace{2\baselineskip}}
\label{fig:holesmesh}
\end{subfigure}
\caption{Domain with boundary conditions and mesh for the asymmetrical bending test case.}
\end{figure}

\begin{figure}
\begin{subfigure}{0.45\linewidth}
\includegraphics[width=\textwidth]{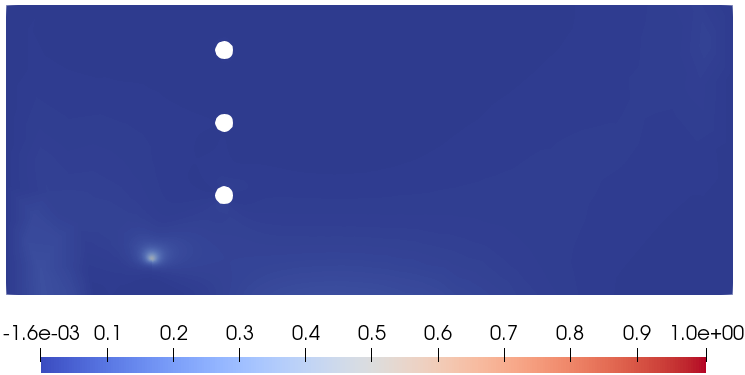}
\caption{{\bf Asymmetrical bending test:} Solution before crack propagation at loading step 11.}
\label{fig:sol-holes-11}
\end{subfigure}
\hspace{0.10\linewidth}
\begin{subfigure}{0.45\linewidth}
\includegraphics[width=\textwidth]{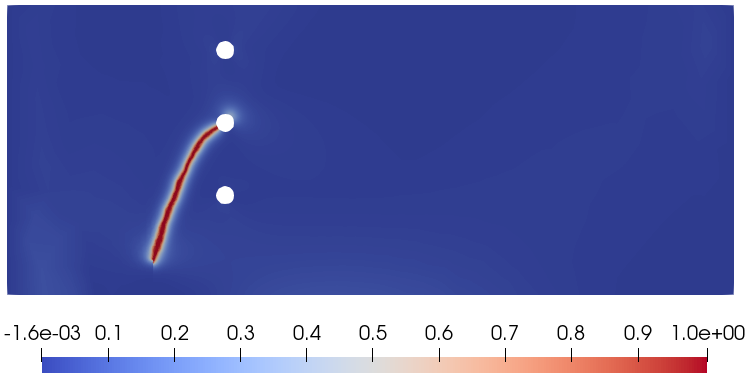}
\caption{{\bf Asymmetrical bending test:} Solution after crack propagation at loading step 12. }
\label{fig:sol-holes-12}
\end{subfigure}
\caption{Solution for $\varphi$ for the asymmetrical bending test case.}
\label{fig:sol-holes}
\end{figure}

This test case considers a rectangular domain with three holes, slightly to the left, and a notch in the lower left part of the domain. It is subject to symmetrical displacement loading on the top boundary,
\begin{equation}
 \label{eq:holesbc}
 \bm u_h^n|_{\Gamma^{Top}} = \begin{pmatrix} 0 \\ \bar{u} n \end{pmatrix}.
\end{equation}
The beam is simply supported as shown in Figure~\ref{fig:holes-domain}. See Figure~\ref{fig:holes-domain} for details on boundary conditions and domain. Experimental results from \cite{holesexperiment} have shown that the crack path should hit the second hole, and we see from the numerical solution, Figure~\ref{fig:sol-holes}, that this also happens here. The mesh has been refined in the region where the crack is expected to propagate, see Figure~\ref{fig:holesmesh}. The problem parameters are chosen similarily to \cite{brun,miehethermodynamics,ambati}, and are presented in Table~\ref{tab:holes}.

\begin{figure}
\begin{subfigure}{0.45\linewidth}
\includegraphics[width = \textwidth]{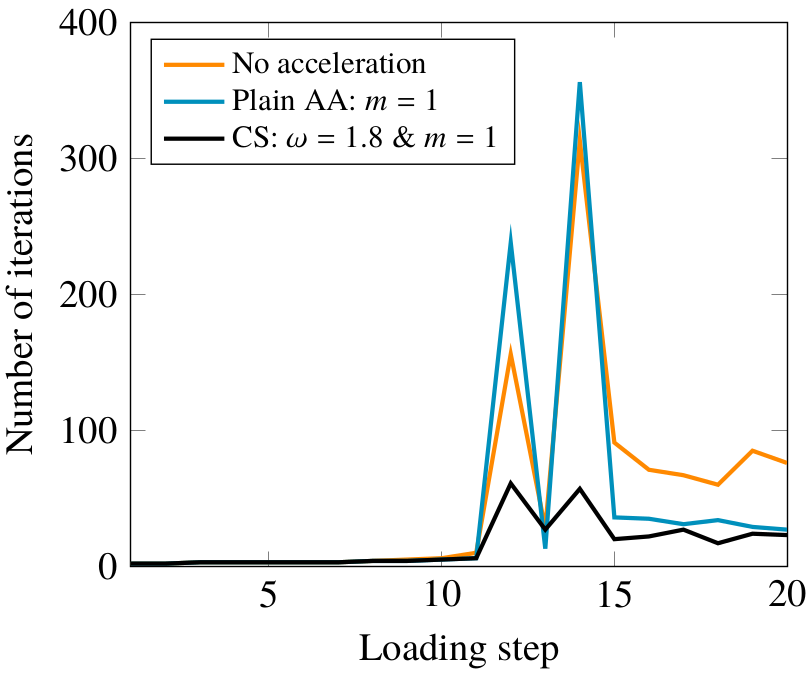}
\caption{Number of iterations per loading step.\\
\vspace{1.3\baselineskip}}
\label{fig:holes-iterations}
\end{subfigure}
\hspace{0.1\linewidth}
\begin{subfigure}{0.45\linewidth}
\includegraphics[width = \textwidth]{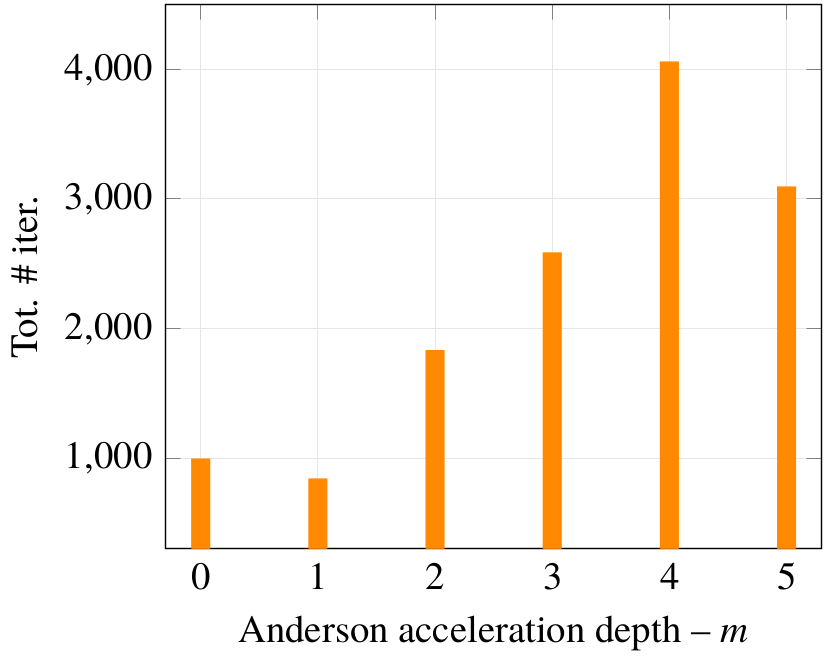}
\caption{Total number of iterations for several depths of plain Anderson acceleration.}
\label{fig:holes-AA}
\end{subfigure}
\caption{{\bf Asymmetrical bending test:} Number of iterations per loading step, and total number of iterations for several depths of Anderson acceleration. ``CS``\ is the combined acceleration scheme. ''Plain AA'' means that Anderson acceleration is applied without any form of safeguard or combination with relaxation.}
\end{figure}

Here, we have two ``critical'' loading steps, in which the crack evolves and a large number of iterations is required, see Figure~\ref{fig:holes-iterations}. For these loading steps, we see that the plain Anderson acceleration does not accelerate, while the combined acceleration performs very well. 

\begin{figure}
\includegraphics[width = \textwidth]{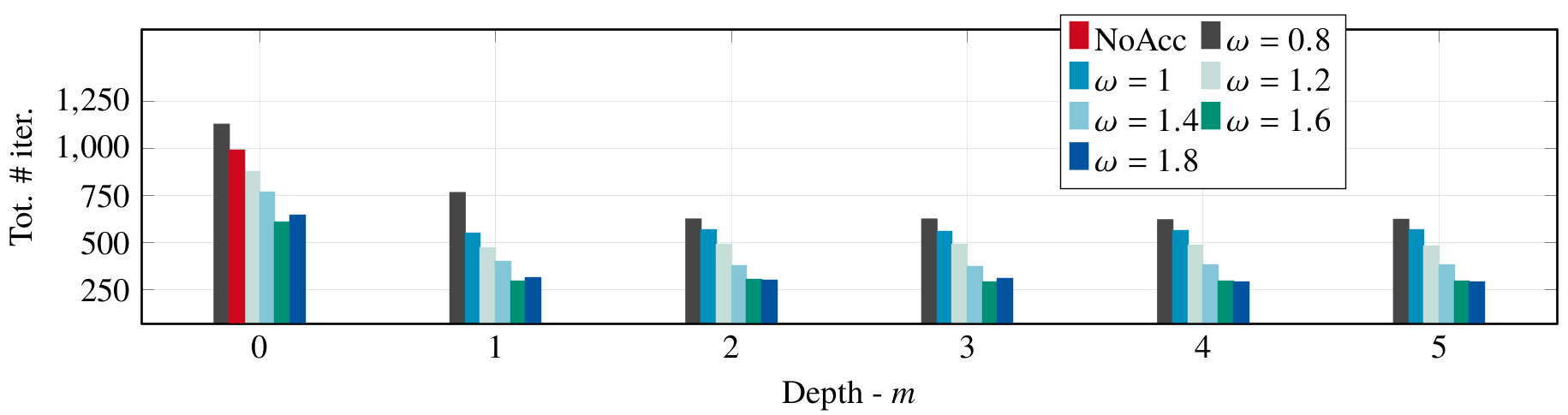}
 \caption{{\bf Asymmetrical bending test:} Total number of iterations for different relaxation parameters and Anderson acceleration depths. ``NoAcc'' is the unaccelerated staggered scheme.}
 \label{fig:holes-totaliterations}
\end{figure}

In Figure~\ref{fig:holes-AA} the total number of iterations for the plain Anderson acceleration is displayed for several depths. We clearly observe that the staggered scheme is significantly decelerated for depths larger than one. In other words, Anderson acceleration is not a robust method in itself for this problem. The combined scheme, on the other hand, reduces the total number of iterations for all combinations of over-relaxation and Anderson acceleration, see Figure~\ref{fig:holes-totaliterations}. There is, however, a tendency that larger relaxation parameters accelerate more, which is expected due to the brutal nature of the crack propagation in loading step 12, see Figure~\ref{fig:sol-holes}.  

The traction vector \eqref{eq:traction} is here calculated on the top boundary and the component of interest is $\tau_y$. In Figure~\ref{fig:holes-load}, the load-displacement curve is displayed, and the displacement is calculated at the left corner of the top boundary. They are, as expected, overlapping as there are no loading steps for any configurations in which the convergence is not achieved in the given maximal amount of iterations. 

\section{Conclusion}\label{sec:conclusion}
The staggered solution scheme is, due to its robustness, a popular method for solving variational phase-field models of brittle fracture. As it often requires a large number of iterations to converge we have proposed a method to accelerate it that exploits the complementary advantages of Anderson acceleration and over-relaxation. The acceleration method alternates between Anderson acceleration and over-relaxation according to a switch that depends on the norms of the previous residuals of the scheme. For problems without brutal crack growth, Anderson acceleration is quite efficient. It is, however, unstable for problems with brutal crack growth, and therefore, not a technique that can be applied without modifications. Over-relaxation, on the other hand, works well within regimes of brutal crack propagation, but might struggle when the iterates get close to the solution. The scheme shows robustness with respect to the tuning parameters, Anderson acceleration depth and relaxation parameter, and converges for all combinations. Moreover, there is a tendency that Anderson acceleration depths larger than one are insignificant, and that over-relaxation with parameters of at least $1.6$ are the best choices. Therefore, we propose to apply the method with depth one and over-relaxation $1.6$, although one might gain some speed in tuning these parameters to specific problems. 

\bibliographystyle{unsrt}
\bibliography{PaperPhaseFieldsArXivNoTemplate.bib}

\end{document}